\documentclass[A4j,11pt]{article}
\usepackage{tikz}
\usepackage{amssymb,amsmath,amsthm,amscd,amsfonts}
\usepackage{amsmath,amssymb,amsthm,amscd,ascmac,mathrsfs,bm,type1cm,multirow,array,enumerate,mathrsfs,tabularx}
\usepackage{graphicx}
\usepackage{graphics}

\newcommand{\ve}{\bm{e}} 

\newcommand{\vv}{{\bm{v}}} 
\newcommand{\vw}{{\bm{w}}} 
\newcommand{\vx}{{\bm{x}}} 

\newcommand{\vZ}{\bm{Z}}
\newcommand{\vP}{\bm{P}}

\begin{document}

\title{{\bf Calabi-Yau structures\\
on the complexifications of\\
rank two symmeric spaces}}
\author{{\bf Naoyuki Koike}}
\date{}
\maketitle

\begin{abstract}
For a (Reimannian) symmetric space $G/K$ of compact type, the natural action of $G$ on its complexification $G^{\mathbb C}/K^{\mathbb C}$ (which is 
an anti-Kaehler symmetric space) is one of the isometric actions called ``Hermann type action''.  
Let $\psi$ be the $G$-invariant strictly plurisubharmonic $C^{\infty}$-function on an open set of $G^{\mathbb C}/K^{\mathbb C}$ arising from 
a $W$-invariant strictly convex $C^{\infty}$-function $\rho$ on an open set of a maximal abelian subspace $\mathfrak a^d$ of $\mathfrak p^d$, 
where $\mathfrak p^d$ is the subspace of the Lie algebra $\mathfrak g^d$ of $G^d$ such that $\mathfrak g^d=\mathfrak k\oplus\mathfrak p^d$ gives 
the Cartan decomposition associated to the dual symmetric space $G^d/K$ of $G/K$ and $W$ is the Weyl group associated to $\mathfrak a^d$.  
In this paper, we first give a new proof of a known relation between the complex Hessian of $\psi$ and the Hessian of $\rho$.  
This new proof is performed from the viewpoint of the orbit geometry of the Hermann type action $G\curvearrowright G^{\mathbb C}/K^{\mathbb C}$ and 
and the isotropy action $K\curvearrowright G^d/K$.  
Next we prove conceptionally that there exists a $C^{\infty}$-Calabi-Yau structure on the whole of the complexification $G^{\mathbb C}/K^{\mathbb C}$ in the case where 
$G/K$ is of rank two on the basis of this relation.  
\end{abstract}


\vspace{0.5truecm}

\section{Introduction}
In this paper, a {\it $C^{\infty}$-Calabi-Yau structure} 
of a $2n$-dimensional manifold $M$ means a quadruple $(J,g,\omega,\Omega)$ 
satisfying the following three conditions:

(i)\ \ $(J,g)$ is a $C^{\infty}$-K$\ddot{\rm a}$hler structure of $M$;

(ii)\ \ $\omega$ is the $C^{\infty}$-K$\ddot{\rm a}$hler form of $(J,g)$;

(iii)\ \ $\Omega$ is a nonvanishing holomorphic $(n,0)$-form on $M$ satisfying 
$$\omega^n=(-1)^{n(n-1)/2}n!\left(\frac{\sqrt{-1}}{2}\right)^n\Omega\wedge\overline{\Omega}\quad({\rm on}\,\,M).\leqno{(1.1)}$$
It follows from $(1.1)$ that $(M,J,g)$ is a $C^{\infty}$-Ricci-flat K$\ddot{\rm a}$hler manifold 
and that, for any real constant $\theta$, a $n$-form ${\rm Re}(e^{\sqrt{-1}\theta}\Omega)$ is a calibration on $(M,g)$, that is, 
${\rm Re}(e^{\sqrt{-1}\theta}\Omega)_p(\ve_1,\cdots,\ve_n)\leq 1$ holds for any $p\in M$ and any orthonormal $n$-system $(\ve_1,\cdots,\ve_n)$ 
of $(T_pM,g_p)$.  Hence, a special Lagrangian submanifold 
$L$ of phase $\theta$ in $(M,J,g,\omega,\Omega)$ is defined as submanifolds calibrated by ${\rm Re}(e^{\sqrt{-1}\theta}\Omega)$, 
that is, ${\rm Re}(e^{\sqrt{-1}\theta}\Omega)_p(\ve_1,\cdots,\ve_n)=1$ holds for any $p\in L$ and any orthonormal frame $(\ve_1,\cdots,\ve_n)$ of 
$(T_pL,g_p\vert_{T_pL\times T_pL}))$.  

The complexification of $C^{\omega}$-Riemannian manifold is defined as follows, where $C^{\omega}$ means the real analyticity.  
Let $(M,g)$ be a $C^{\omega}$-Riemannian manifold, $TM$ be the tangent bundle of $M$ and $J_A$ be the adapted complex structure (associated to $g$) defined on 
a tubular neighborhood $U_A$ of the zero section of $TM$, where we note that $J_A$ is defined on the whole of $TM$ (i.e., $U_A=TM$) when $(M,g)$ is of non-negative 
curvature.  See \cite{Sz} and \cite{GS1,GS2} about the definition of the adapted complex structure.  
The complex manifold $(U_A,J_A)$ is regarded as the complexification of $(M,g)$ under the identification of 
the zero-section of $TM$ with $M$.  Let $G/K$ be a symmetric space of compact type.  
Since $G/K$ is a $C^{\omega}$-Riemannian manifold of non-negative curvature, 
the adapted complex structure $J_A$ of $T(G/K)$ is defined on the whole of the tangent bundle $T(G/K)$.  
As above, $(T(G/K),J_A)$ is regarded as the complexification of $G/K$.  
We can also define the complexification of $G/K$ as follows.  

Let $G^{\mathbb C}$ and $K^{\mathbb C}$ be the complexifications of $G$ and $K$, respectively.  
Denote by $\mathfrak g,\mathfrak k,\mathfrak g^{\mathbb C}$ and $\mathfrak k^{\mathbb C}$ be the Lie algebras of 
$G,K,G^{\mathbb C}$ and $K^{\mathbb C}$, respectively.  
Let $\mathfrak g=\mathfrak k\oplus\mathfrak p$ and 
$\mathfrak g^{\mathbb C}=\mathfrak k^{\mathbb C}\oplus\mathfrak p^{\mathbb C}$ be the canonical decompositions 
associated to the semi-simple symmetric pairs $(G,K)$ and $(G^{\mathbb C},K^{\mathbb C})$, respectively.  
Here we note that $\mathfrak p$ and $\mathfrak p^{\mathbb C}$ are identified with the tangent spaces 
$T_{eK}(G/K)$ and $T_{eK^{\mathbb C}}(G^{\mathbb C}/K^{\mathbb C})$, respectively, where $e$ is the identity element of $G^{\mathbb C}$.  
For the simplicity, set $o:=eK^{\mathbb C}(=eK)$.  
Define a complex linear transformation ${\bf j}$ of $\mathfrak p^{\mathbb C}$ by 
${\bf j}(v):=\sqrt{-1}v\,\,\,\,(v\in\mathfrak p^{\mathbb C})$.  Since ${\bf j}$ is 
${\rm Ad}_{G^{\mathbb C}}(K^{\mathbb C})$-invariant, we can define $G^{\mathbb C}$-invariant complex structure 
${\bf J}$ of $G^{\mathbb C}/K^{\mathbb C}$ satisfying ${\bf J}_o={\bf j}$ uniquely, where 
${\rm Ad}_{G^{\mathbb C}}$ denotes the adjoint representation of $G^{\mathbb C}$.  
The complex manifold $(G^{\mathbb C}/K^{\mathbb C},{\bf J})$ is regarded as another complexification of $G/K$ 
under the identification of the orbit $G\cdot o$ of the subaction $G\curvearrowright G^{\mathbb C}/K^{\mathbb C}$ of the natural action 
$G^{\mathbb C}\curvearrowright G^{\mathbb C}/K^{\mathbb C}$ with $G/K$.  
The subaction $G\curvearrowright G^{\mathbb C}/K^{\mathbb C}$ is one of the actions called a {\it Hermann type action} (in \cite{K2}).  
A Hermann type action is defined as follows.  Let $H_1$ and $H_2$ be symmetric subgroups of a semi-simple Lie group $G$, that is, closed subgroups of $G$ 
admitting involutions $\sigma_i$ of $G$ with $({\rm Fix}\,\sigma_i)_0\subset H_i\subset{\rm Fix}\,\sigma_i$, where ${\rm Fix}\,\sigma_i$ is the fixed point group of 
$\sigma_i$ and $({\rm Fix}\,\sigma_i)_0$ is the identity component of ${\rm Fix}\,\sigma_i$.  
In particular, if $G$ is compact, then this action $H_2\curvearrowright G/H_1$ is called a {\it Hermann action}.  
The geometry of the orbits of Hermann actions was investigated by O. Goetsches and G. Thorbergsson (\cite{GT}) and 
the geometry of the orbits of Hermann type actions was investigated by the author (\cite{K2}).  
We can define a natural holomorphic diffeomorphism between two complexifications $(T(G/K),J_A)$ and 
$(G^{\mathbb C}/K^{\mathbb C},{\bf J})$ of $G/K$ as follows.  
Let $B$ be the Killing form of $\mathfrak g$ and set $B_A:=-2{\rm Re}\,B^{\mathbb C}$.  
Since $B_A$ is ${\rm Ad}_{G^{\mathbb C}}(K^{\mathbb C})$-invariant, we can define $G^{\mathbb C}$-invariant 
pseudo-Riemannian metric $\beta_A$ of $G^{\mathbb C}/K^{\mathbb C}$ satisfying $(\beta_A)_o=B_A$ uniquely.  
The pseudo-Riemannian manifold $(G^{\mathbb C}/K^{\mathbb C},\beta_A)$ is one of semi-simple 
pseudo-Riemannian symmetric spaces.  Also, the triple $(G^{\mathbb C}/K^{\mathbb C},{\bf J},\beta_A)$ gives 
an anti-K$\ddot{\rm a}$hler manifold.  See \cite{BFV} (\cite{K1} also) about the definition of the anti-K$\ddot{\rm a}$hler manifold.  
Denote by ${\rm Exp}_p$ the exponential map of $(G^{\mathbb C}/K^{\mathbb C},\beta_A)$ at 
$p\in G^{\mathbb C}/K^{\mathbb C}$.  The natural bijection 
$\displaystyle{\Psi:(T(G/K),J_A)\mathop{\longrightarrow}_{\cong}(G^{\mathbb C}/K^{\mathbb C},{\bf J})}$ is given by 
$$\Psi(v):={\rm Exp}_p({\bf J}_p(v))\quad\,\,(p\in G/K,\,\,v\in T_p(G/K))$$
(see Figure 1), where $v$ in the right-hand side is regarded as a tangent vector of the submanifold $G\cdot o(\approx G/K)$ in 
$G^{\mathbb C}/K^{\mathbb C}$.  
It is shown that $\Psi$ is a holomorphic diffeomorphism between $(T(G/K),J_A)$ and 
$(G^{\mathbb C}/K^{\mathbb C},{\bf J})$ (see Section 2 about the proof of this fact).  
Thus these two complexifications of $G/K$ are identified through $\Psi$.  

\vspace{0.5truecm}

\centerline{
\unitlength 0.1in
\begin{picture}( 32.1400, 17.2900)(  8.1000,-23.5900)
%
\special{pn 8}%
\special{ar 3300 1466 382 94  6.2831853 6.2831853}%
\special{ar 3300 1466 382 94  0.0000000 3.1415927}%
%
\special{pn 8}%
\special{ar 3300 1466 382 94  3.1415927 3.1922256}%
\special{ar 3300 1466 382 94  3.3441243 3.3947572}%
\special{ar 3300 1466 382 94  3.5466559 3.5972889}%
\special{ar 3300 1466 382 94  3.7491876 3.7998205}%
\special{ar 3300 1466 382 94  3.9517192 4.0023521}%
\special{ar 3300 1466 382 94  4.1542509 4.2048838}%
\special{ar 3300 1466 382 94  4.3567825 4.4074154}%
\special{ar 3300 1466 382 94  4.5593142 4.6099471}%
\special{ar 3300 1466 382 94  4.7618458 4.8124787}%
\special{ar 3300 1466 382 94  4.9643775 5.0150104}%
\special{ar 3300 1466 382 94  5.1669091 5.2175420}%
\special{ar 3300 1466 382 94  5.3694408 5.4200737}%
\special{ar 3300 1466 382 94  5.5719724 5.6226053}%
\special{ar 3300 1466 382 94  5.7745040 5.8251370}%
\special{ar 3300 1466 382 94  5.9770357 6.0276686}%
\special{ar 3300 1466 382 94  6.1795673 6.2302002}%
%
\special{pn 8}%
\special{ar 1908 1466 1014 806  5.5340129 6.2831853}%
\special{ar 1908 1466 1014 806  0.0000000 0.7853982}%
%
\special{pn 8}%
\special{ar 4694 1466 1012 806  2.3561945 3.8903018}%
%
\special{pn 8}%
\special{ar 4346 1552 786 942  2.5896390 3.8667463}%
%
\special{pn 8}%
\special{ar 2288 1552 784 942  5.5599686 6.2831853}%
\special{ar 2288 1552 784 942  0.0000000 0.5503859}%
%
\special{pn 20}%
\special{sh 1}%
\special{ar 3560 1534 10 10 0  6.28318530717959E+0000}%
\special{sh 1}%
\special{ar 3560 1534 10 10 0  6.28318530717959E+0000}%
%
\special{pn 20}%
\special{sh 1}%
\special{ar 3066 1538 10 10 0  6.28318530717959E+0000}%
\special{sh 1}%
\special{ar 3066 1538 10 10 0  6.28318530717959E+0000}%
%
\special{pn 8}%
\special{pa 3952 1080}%
\special{pa 3732 988}%
\special{dt 0.045}%
\special{sh 1}%
\special{pa 3732 988}%
\special{pa 3786 1032}%
\special{pa 3780 1008}%
\special{pa 3800 994}%
\special{pa 3732 988}%
\special{fp}%
\put(30.7900,-15.2100){\makebox(0,0)[lb]{$o$}}%
\put(35.8600,-15.6300){\makebox(0,0)[lt]{$p$}}%
\put(40.2400,-17.1500){\makebox(0,0)[lt]{{\small $G\cdot o=G/K$}}}%
\put(39.8500,-10.7500){\makebox(0,0)[lt]{{\small $\Psi(T_p(G\cdot o))$}}}%
\put(35.0000,-13.6000){\makebox(0,0)[rb]{${\bf J}_p(v)$}}%
\put(38.0100,-14.4000){\makebox(0,0)[lt]{$v$}}%
\put(31.5500,-21.4000){\makebox(0,0)[lt]{$G^{\mathbb C}/K^{\mathbb C}$}}%
\put(28.9700,-23.5900){\makebox(0,0)[lt]{$\,$}}%
%
\special{pn 13}%
\special{pa 3568 1538}%
\special{pa 3782 1478}%
\special{fp}%
\special{sh 1}%
\special{pa 3782 1478}%
\special{pa 3712 1476}%
\special{pa 3732 1492}%
\special{pa 3724 1514}%
\special{pa 3782 1478}%
\special{fp}%
%
\special{pn 13}%
\special{pa 3554 1538}%
\special{pa 3560 1272}%
\special{fp}%
\special{sh 1}%
\special{pa 3560 1272}%
\special{pa 3540 1338}%
\special{pa 3560 1326}%
\special{pa 3578 1340}%
\special{pa 3560 1272}%
\special{fp}%
%
\special{pn 20}%
\special{sh 1}%
\special{ar 3606 1242 10 10 0  6.28318530717959E+0000}%
\special{sh 1}%
\special{ar 3606 1242 10 10 0  6.28318530717959E+0000}%
%
\special{pn 8}%
\special{pa 3858 1310}%
\special{pa 3612 1242}%
\special{dt 0.045}%
\special{sh 1}%
\special{pa 3612 1242}%
\special{pa 3670 1278}%
\special{pa 3662 1256}%
\special{pa 3682 1240}%
\special{pa 3612 1242}%
\special{fp}%
\put(38.8300,-12.9200){\makebox(0,0)[lt]{{\small $\Psi(v)={\rm Exp}_p({\bf J}_p(v))$}}}%
%
\special{pn 8}%
\special{ar 4044 1566 668 206  1.6326912 1.7699909}%
\special{ar 4044 1566 668 206  1.8523708 1.9896706}%
\special{ar 4044 1566 668 206  2.0720504 2.2093502}%
\special{ar 4044 1566 668 206  2.2917301 2.4290298}%
\special{ar 4044 1566 668 206  2.5114097 2.6487095}%
\special{ar 4044 1566 668 206  2.7310893 2.8683891}%
\special{ar 4044 1566 668 206  2.9507690 3.0203158}%
%
\special{pn 8}%
\special{pa 3384 1596}%
\special{pa 3370 1566}%
\special{da 0.070}%
\special{sh 1}%
\special{pa 3370 1566}%
\special{pa 3382 1636}%
\special{pa 3394 1614}%
\special{pa 3418 1618}%
\special{pa 3370 1566}%
\special{fp}%
%
\special{pn 8}%
\special{ar 1646 1514 382 94  6.2831853 6.2831853}%
\special{ar 1646 1514 382 94  0.0000000 3.1415927}%
%
\special{pn 8}%
\special{ar 1646 1514 382 94  3.1415927 3.1922256}%
\special{ar 1646 1514 382 94  3.3441243 3.3947572}%
\special{ar 1646 1514 382 94  3.5466559 3.5972889}%
\special{ar 1646 1514 382 94  3.7491876 3.7998205}%
\special{ar 1646 1514 382 94  3.9517192 4.0023521}%
\special{ar 1646 1514 382 94  4.1542509 4.2048838}%
\special{ar 1646 1514 382 94  4.3567825 4.4074154}%
\special{ar 1646 1514 382 94  4.5593142 4.6099471}%
\special{ar 1646 1514 382 94  4.7618458 4.8124787}%
\special{ar 1646 1514 382 94  4.9643775 5.0150104}%
\special{ar 1646 1514 382 94  5.1669091 5.2175420}%
\special{ar 1646 1514 382 94  5.3694408 5.4200737}%
\special{ar 1646 1514 382 94  5.5719724 5.6226053}%
\special{ar 1646 1514 382 94  5.7745040 5.8251370}%
\special{ar 1646 1514 382 94  5.9770357 6.0276686}%
\special{ar 1646 1514 382 94  6.1795673 6.2302002}%
%
\special{pn 8}%
\special{pa 1266 926}%
\special{pa 1266 2048}%
\special{fp}%
%
\special{pn 8}%
\special{pa 2028 926}%
\special{pa 2028 2048}%
\special{fp}%
%
\special{pn 8}%
\special{pa 2380 1500}%
\special{pa 2638 1500}%
\special{fp}%
\special{sh 1}%
\special{pa 2638 1500}%
\special{pa 2570 1480}%
\special{pa 2584 1500}%
\special{pa 2570 1520}%
\special{pa 2638 1500}%
\special{fp}%
\put(24.6000,-14.6000){\makebox(0,0)[lb]{$\Psi$}}%
%
\special{pn 8}%
\special{ar 1174 1682 440 214  6.2831853 6.4669525}%
\special{ar 1174 1682 440 214  6.5772129 6.7609801}%
\special{ar 1174 1682 440 214  6.8712404 7.0550077}%
\special{ar 1174 1682 440 214  7.1652680 7.3490352}%
\special{ar 1174 1682 440 214  7.4592956 7.6430628}%
\special{ar 1174 1682 440 214  7.7533231 7.8539816}%
%
\special{pn 8}%
\special{pa 1622 1682}%
\special{pa 1638 1612}%
\special{da 0.070}%
\special{sh 1}%
\special{pa 1638 1612}%
\special{pa 1604 1672}%
\special{pa 1626 1664}%
\special{pa 1642 1680}%
\special{pa 1638 1612}%
\special{fp}%
\put(8.1000,-23.5900){\makebox(0,0)[lt]{{\small $0$-section(=G/K)}}}%
\put(14.4700,-21.4500){\makebox(0,0)[lt]{$T(G/K)$}}%
%
\special{pn 20}%
\special{sh 1}%
\special{ar 1878 1594 10 10 0  6.28318530717959E+0000}%
\special{sh 1}%
\special{ar 1878 1594 10 10 0  6.28318530717959E+0000}%
%
\special{pn 13}%
\special{pa 1870 1584}%
\special{pa 1876 1318}%
\special{fp}%
\special{sh 1}%
\special{pa 1876 1318}%
\special{pa 1854 1384}%
\special{pa 1874 1372}%
\special{pa 1894 1386}%
\special{pa 1876 1318}%
\special{fp}%
\put(18.2800,-13.9700){\makebox(0,0)[rb]{$v$}}%
%
\special{pn 8}%
\special{pa 1878 988}%
\special{pa 1878 2110}%
\special{fp}%
\put(14.4700,-8.0000){\makebox(0,0)[lb]{$T_p(G/K)$}}%
\put(18.4400,-16.2900){\makebox(0,0)[rt]{$p$}}%
%
\special{pn 8}%
\special{ar 1174 2288 224 400  3.0696708 3.2619785}%
\special{ar 1174 2288 224 400  3.3773631 3.5696708}%
\special{ar 1174 2288 224 400  3.6850554 3.8773631}%
\special{ar 1174 2288 224 400  3.9927477 4.1850554}%
\special{ar 1174 2288 224 400  4.3004400 4.4927477}%
\special{ar 1174 2288 224 400  4.6081324 4.7123890}%
%
\special{pn 8}%
\special{ar 1878 846 274 250  1.7019292 1.9318143}%
\special{ar 1878 846 274 250  2.0697453 2.2996304}%
\special{ar 1878 846 274 250  2.4375614 2.6674465}%
\special{ar 1878 846 274 250  2.8053775 3.0352626}%
%
\special{pn 8}%
\special{pa 1836 1086}%
\special{pa 1878 1086}%
\special{fp}%
\special{sh 1}%
\special{pa 1878 1086}%
\special{pa 1812 1066}%
\special{pa 1826 1086}%
\special{pa 1812 1106}%
\special{pa 1878 1086}%
\special{fp}%
\end{picture}%
\hspace{1.5truecm}}

\vspace{0.5truecm}

\centerline{{\bf Figure 1$\,$:$\,$ The identification of $(T(G/K),J_A)$ and $(G^{\mathbb C}/K^{\mathbb C},{\bf J})$}}

\vspace{0.5truecm}

Let $\psi$ be the $G$-invariant strictly plurisubharmonic $C^{\infty}$-function on an open set of $G^{\mathbb C}/K^{\mathbb C}$ arising from 
a $W$-invariant strictly convex $C^{\infty}$-function $\rho$ on an open set of a maximal abelian subspace $\mathfrak a^d$ of $\mathfrak p^d$, 
where $\mathfrak p^d$ is the subspace of the Lie algebra $\mathfrak g^d$ of $G^d$ such that $\mathfrak g^d=\mathfrak k\oplus\mathfrak p^d$ gives 
the Cartan decomposition associated to the dual symmetric space $G^d/K$ of $G/K$ and $W$ is the Weyl group associated to $\mathfrak a^d$.  
T. Delcroix (\cite{D1}) derived relations between the complex Hessian of $\psi$ and the Hessian of $\rho$ 
in the case where $G^{\mathbb C}/K^{\mathbb C}$ is the complexification of a compact semi-simple Lie group.  
Later, T. Delcroix (\cite{D2}) showed that the generalizations of the relations hold for horosymmetric homogeneous spaces 
(see the relations in Corollaries 3.11 and 3.13 of \cite{D2}).
In this paper, we give a new proof of the relations (see (3,17) and $(3.18)$ in Lemma 3.1) in the case where $G^{\mathbb C}/K^{\mathbb C}$ is the complexification 
of a general symmetric space $G/K$ of compact type.  This new proof is performed from the viewpoint of the orbit geometry of the Hermann type action 
$G\curvearrowright G^{\mathbb C}/K^{\mathbb C}$.  In more detail, it is performed by using the explicit descriptions of the shape operators of the orbits of 
the isotropy action $K\curvearrowright G^d/K$ and the Hermann type action $G\curvearrowright G^{\mathbb C}/K^{\mathbb C}$ (see Section 3).  

In 1993, M. B. Stenzel (\cite{St}) showed that there exists a $G$-invariant complete Ricci-flat metric on the whole of the cotangent bundle of a rank one symmetric space $G/K$ 
of compact type.  
In 2004, R. Bielawski (\cite{B3}) investigated the existence of $G$-invariant metrics with the prescribed Ricci tensor on the complexification 
$G^{\mathbb C}/K^{\mathbb C}$ of a general rank symmetric space $G/K$ of compact type (see \cite{B1} also), where $G,\,K$ and $G^{\mathbb C}$ are assumed to be connected.  
Note that, in \cite{B3}, it is stated that the $G$-invariant metrics with the prescribed Ricci tensor are defined on the {\it whole of} $G^{\mathbb C}/K^{\mathbb C}$ 
but it seems that there is a gap in the proof of the fact that the $G$-invariant metrics are defined on the {\it whole of} $G^{\mathbb C}/K^{\mathbb C}$ 
(It seems that there is a gap in the process of the proof of the main theorem in \cite{B1} also).  
In 2019, the author (\cite{K4}) investigated the existence of Calabi-Yau structures on the complexification $G^{\mathbb C}/K^{\mathbb C}$ of a general rank symmetric space 
$G/K$ of compact type 
but there are some gaps in the proof of some lemmas in the paper.  
In this paper, we close these gaps and prove the following fact in the case where $G/K$ is of rank two.  

\vspace{0.5truecm}

\noindent
{\bf Theorem A.} {\sl If $G/K$ is a rank two symmetric space of compact type, then there exists a $G$-invariant $C^{\infty}$-Calabi-Yau structure 
on the whole of the complexification $G^{\mathbb C}/K^{\mathbb C}$ of $G/K$.}

\vspace{0.5truecm}

\noindent
{\it Remark 1.1.}\ \ In the case where the restricted root system of the symmetric pair $(G,K)$ is other than $(\mathfrak g_2)$-type, 
T.T. Nghiem (\cite{N}) constructed $G$-invariant complete Calabi-Yau metrics on the whole of $G^{\mathbb C}/K^{\mathbb C}$ of Euclidean volume growth 
having any prescribed singular $G^{\mathbb C}$-horospherical tangent cone as its tangent cone at the infinity.  
By examples of $G$-inavriant complete Calabi-Yau metrics constructed by them, it is shown that there exist a $G$-invariant $C^{\infty}$-Calabi-Yau structure 
on the whole of the complexification $G^{\mathbb C}/K^{\mathbb C}$ of $G/K$ in the case where the restricted root system of the symmetric pair $(G,K)$ is other than 
$(\mathfrak g_2)$-type.  However its proof is not conceptional.  On the other hand, our proof is conceptional.  

\vspace{0.5truecm}

The Hermann type action $G\curvearrowright G^{\mathbb C}/K^{\mathbb C}$ is a complex hyperpolar action 
(see \cite{K1} about the definition of the complex hyperpolar action).  So, the following question arises naturally.  

\vspace{0.5truecm}

\noindent
{\bf Question.} {\sl 
Let $G\curvearrowright(M,{\bf J},\beta_A,\Omega)$ be a complex hyperpolar action of a compact semi-simple Lie group $G$ 
on an anti-Kaehler manifold $(M,{\bf J},\beta_A,\Omega)$ equipped with a nowhere zero holomorphic $(n,0)$-form $\Omega$.  Assume that 
a $G$-orbit $G\cdot o$ gives a real form of $(M,{\bf J},\beta_A,\Omega)$.  
In what case, can we show the existence of a $G$-invariant Riemannian metric $\beta$ on $(M,{\bf J},\Omega)$ such that $({\bf J},\beta,\omega,\Omega)$ 
is a Calabi-Yau structure?  Here $\omega$ denotes the symplectic form given by ${\bf J}$ and $\beta$.  
}  

\section{Basic notions and facts} 
In this section, we recall some basic notions and facts.  

\vspace{0.25truecm}

\noindent
{\bf 2.1.\ Adapted complex structure}\ \ 
We recall the adapted complex structure of the tangent bundle of complete 
$C^{\omega}$-Riemanninan manifold $(M,g)$, which was introduced by R. Sz$\ddot{\rm o}$ke (\cite{Sz}) and 
V. Guillemin-M. B. Stenzel (\cite{GS1,GS2}), 
and show that $\Psi$ stated in Introduction is a holomorphic diffeomorphism.  
For each geodesic $\gamma:\mathbb R\to M$ in $(M,g)$, we define $\gamma_{\ast}:\mathbb C\to TM$ by 
$$\gamma_{\ast}(s+t\sqrt{-1}):=t\gamma'(s)\quad\,\,(s+t\sqrt{-1}\in\mathbb C),$$
where $TM$ is the tangent bunndle of $M$ and $\gamma'(s)$ is the velocity vector of $\gamma$ at $s$.  
Then V. Guillemin and M. B. Stenzel showed that 
there exists a unique complex structure $J_A$ of $TM$ such that, for any geodesic $\gamma:\mathbb R\to M$ in 
$(M,g)$, $\gamma_{\ast}:\mathbb C\to(TM,J_A)$ is holomorphic.  This complex structure is called the 
{\it adapted complex structure} of $TM$ associated to $g$.  

\vspace{0.5truecm}

\noindent
{\bf Proposition 2.1.} {\sl Let $G/K$ be a symmetric space of compact type and 
$\Psi:T(G/K)\to G^{\mathbb C}/K^{\mathbb C}$ be the map stated in Introduction.  
Then the map $\Psi$ is a holomorphic diffeomorphism of $(T(G/K),J_A)$ onto $(G^{\mathbb C}/K^{\mathbb C},{\bf J})$.}

\vspace{0.5truecm}

\noindent
{\it Proof.} It is clear that $\Psi$ is a diffeomorphism.  We shall show that $\Psi$ is holomorphic.  
Let $\gamma:\mathbb R\to G/K$ be a geodesic in $G/K$ and set $\gamma^{\mathbb C}:=\Psi\circ\gamma_{\ast}$.  
Then we have 
$$\gamma^{\mathbb C}(s+t\sqrt{-1})={\rm Exp}_{\gamma(s)}(t{\bf J}_{\gamma(s)}(\gamma'(s)))$$
(see Figure 2), which is the complexification of the geodesic $\gamma$ in 
$G/K(\approx G\cdot o\subset G^{\mathbb C}/K^{\mathbb C}))$ in the sense of \cite{K2,K3}.  
By using Proposition 3.1 of \cite{K3}, we can show that 
$\gamma^{\mathbb C}:\mathbb C\to(G^{\mathbb C}/K^{\mathbb C},{\bf J})$ is holomorphic.  
This fact together with the arbitrariness of $\gamma$ implies that $\Psi$ is holomorphic.  \qed

\vspace{0.35truecm}

\centerline{
\unitlength 0.1in
\begin{picture}( 54.7200, 22.3700)( -5.0000,-23.3500)
%
\special{pn 8}%
\special{ar 2730 1464 584 164  0.0000000 6.2831853}%
%
\special{pn 8}%
\special{pa 2506 1620}%
\special{pa 3020 1326}%
\special{pa 3020 676}%
\special{pa 2506 928}%
\special{pa 2506 928}%
\special{pa 2506 1620}%
\special{fp}%
%
\special{pn 8}%
\special{pa 3020 1326}%
\special{pa 3020 1578}%
\special{fp}%
%
\special{pn 8}%
\special{pa 3020 1636}%
\special{pa 3020 1928}%
\special{fp}%
%
\special{pn 8}%
\special{pa 2506 2222}%
\special{pa 2506 1612}%
\special{fp}%
%
\special{pn 20}%
\special{sh 1}%
\special{ar 2758 1228 10 10 0  6.28318530717959E+0000}%
\special{sh 1}%
\special{ar 2758 1228 10 10 0  6.28318530717959E+0000}%
%
\special{pn 20}%
\special{sh 1}%
\special{ar 2758 984 10 10 0  6.28318530717959E+0000}%
\special{sh 1}%
\special{ar 2758 984 10 10 0  6.28318530717959E+0000}%
%
\special{pn 8}%
\special{pa 1050 968}%
\special{pa 1836 968}%
\special{fp}%
\special{sh 1}%
\special{pa 1836 968}%
\special{pa 1770 948}%
\special{pa 1784 968}%
\special{pa 1770 988}%
\special{pa 1836 968}%
\special{fp}%
%
\special{pn 8}%
\special{pa 1438 1506}%
\special{pa 1438 366}%
\special{fp}%
\special{sh 1}%
\special{pa 1438 366}%
\special{pa 1418 434}%
\special{pa 1438 420}%
\special{pa 1458 434}%
\special{pa 1438 366}%
\special{fp}%
%
\special{pn 8}%
\special{pa 1050 724}%
\special{pa 1836 724}%
\special{fp}%
%
\special{pn 8}%
\special{pa 1040 480}%
\special{pa 1826 480}%
\special{fp}%
%
\special{pn 8}%
\special{pa 1050 1212}%
\special{pa 1836 1212}%
\special{fp}%
%
\special{pn 8}%
\special{pa 1040 1456}%
\special{pa 1826 1456}%
\special{fp}%
%
\special{pn 8}%
\special{pa 2506 1368}%
\special{pa 3020 1090}%
\special{fp}%
%
\special{pn 8}%
\special{pa 2506 1114}%
\special{pa 3020 838}%
\special{fp}%
%
\special{pn 8}%
\special{pa 2506 1838}%
\special{pa 3020 1562}%
\special{fp}%
%
\special{pn 13}%
\special{pa 2758 1482}%
\special{pa 2750 1236}%
\special{fp}%
\special{sh 1}%
\special{pa 2750 1236}%
\special{pa 2732 1304}%
\special{pa 2752 1290}%
\special{pa 2772 1302}%
\special{pa 2750 1236}%
\special{fp}%
%
\special{pn 20}%
\special{sh 1}%
\special{ar 1438 724 10 10 0  6.28318530717959E+0000}%
\special{sh 1}%
\special{ar 1438 480 10 10 0  6.28318530717959E+0000}%
\special{sh 1}%
\special{ar 1438 480 10 10 0  6.28318530717959E+0000}%
%
\special{pn 8}%
\special{ar 4282 1464 584 164  0.0000000 6.2831853}%
%
\special{pn 8}%
\special{pa 4060 1620}%
\special{pa 4574 1326}%
\special{pa 4574 676}%
\special{pa 4060 928}%
\special{pa 4060 928}%
\special{pa 4060 1620}%
\special{fp}%
%
\special{pn 8}%
\special{pa 4574 1326}%
\special{pa 4574 1578}%
\special{fp}%
%
\special{pn 8}%
\special{pa 4574 1636}%
\special{pa 4574 1944}%
\special{fp}%
%
\special{pn 8}%
\special{pa 4060 2222}%
\special{pa 4060 1612}%
\special{fp}%
%
\special{pn 20}%
\special{sh 1}%
\special{ar 4312 1228 10 10 0  6.28318530717959E+0000}%
\special{sh 1}%
\special{ar 4312 1228 10 10 0  6.28318530717959E+0000}%
%
\special{pn 20}%
\special{sh 1}%
\special{ar 4312 984 10 10 0  6.28318530717959E+0000}%
\special{sh 1}%
\special{ar 4312 984 10 10 0  6.28318530717959E+0000}%
%
\special{pn 8}%
\special{pa 4060 1368}%
\special{pa 4574 1090}%
\special{fp}%
%
\special{pn 8}%
\special{pa 4060 1114}%
\special{pa 4574 838}%
\special{fp}%
%
\special{pn 8}%
\special{pa 4060 1838}%
\special{pa 4574 1562}%
\special{fp}%
%
\special{pn 8}%
\special{pa 4060 2058}%
\special{pa 4574 1782}%
\special{fp}%
%
\special{pn 13}%
\special{pa 4312 1482}%
\special{pa 4302 1236}%
\special{fp}%
\special{sh 1}%
\special{pa 4302 1236}%
\special{pa 4286 1304}%
\special{pa 4304 1290}%
\special{pa 4326 1302}%
\special{pa 4302 1236}%
\special{fp}%
%
\special{pn 8}%
\special{pa 3332 1270}%
\special{pa 3768 1270}%
\special{fp}%
\special{sh 1}%
\special{pa 3768 1270}%
\special{pa 3702 1250}%
\special{pa 3716 1270}%
\special{pa 3702 1290}%
\special{pa 3768 1270}%
\special{fp}%
%
\special{pn 8}%
\special{pa 2022 854}%
\special{pa 2322 1042}%
\special{fp}%
\special{sh 1}%
\special{pa 2322 1042}%
\special{pa 2276 990}%
\special{pa 2276 1014}%
\special{pa 2254 1024}%
\special{pa 2322 1042}%
\special{fp}%
\put(35.0600,-12.2800){\makebox(0,0)[lb]{$\Psi$}}%
\put(13.4100,-16.6800){\makebox(0,0)[lt]{$\mathbb C$}}%
\put(25.4500,-23.3500){\makebox(0,0)[lt]{$T(G/K)$}}%
\put(40.7000,-23.3000){\makebox(0,0)[lt]{$G^{\mathbb C}/K^{\mathbb C}$}}%
%
\special{pn 20}%
\special{sh 1}%
\special{ar 1730 968 10 10 0  6.28318530717959E+0000}%
\special{sh 1}%
\special{ar 1730 968 10 10 0  6.28318530717959E+0000}%
%
\special{pn 20}%
\special{sh 1}%
\special{ar 2952 1368 10 10 0  6.28318530717959E+0000}%
\special{sh 1}%
\special{ar 2952 1368 10 10 0  6.28318530717959E+0000}%
%
\special{pn 13}%
\special{pa 4312 1482}%
\special{pa 4496 1376}%
\special{fp}%
\special{sh 1}%
\special{pa 4496 1376}%
\special{pa 4428 1392}%
\special{pa 4450 1402}%
\special{pa 4448 1426}%
\special{pa 4496 1376}%
\special{fp}%
\put(16.9000,-10.0900){\makebox(0,0)[lt]{{\small $1$}}}%
\put(13.9000,-4.4000){\makebox(0,0)[rb]{{\small $2\sqrt{-1}$}}}%
\put(14.0000,-6.9000){\makebox(0,0)[rb]{{\small $\sqrt{-1}$}}}%
\put(18.7500,-10.0100){\makebox(0,0)[lb]{$s$}}%
\put(14.9600,-3.6600){\makebox(0,0)[lb]{$t\sqrt{-1}$}}%
\put(21.5600,-8.8700){\makebox(0,0)[lb]{$\gamma_{\ast}$}}%
\put(27.0000,-12.2000){\makebox(0,0)[rb]{{\small $v$}}}%
%
\special{pn 8}%
\special{ar 3118 1164 174 220  3.3801146 3.6854581}%
\special{ar 3118 1164 174 220  3.8686642 4.1740077}%
\special{ar 3118 1164 174 220  4.3572138 4.6625573}%
%
\special{pn 8}%
\special{pa 2952 1114}%
\special{pa 2952 1132}%
\special{da 0.070}%
\special{sh 1}%
\special{pa 2952 1132}%
\special{pa 2972 1064}%
\special{pa 2952 1078}%
\special{pa 2932 1064}%
\special{pa 2952 1132}%
\special{fp}%
\put(31.5600,-9.9200){\makebox(0,0)[lb]{{\small $\gamma'(s)$}}}%
%
\special{pn 8}%
\special{ar 3118 878 174 220  3.3801146 3.6854581}%
\special{ar 3118 878 174 220  3.8686642 4.1740077}%
\special{ar 3118 878 174 220  4.3572138 4.6625573}%
\put(31.4700,-7.1600){\makebox(0,0)[lb]{{\small $2\gamma'(s)$}}}%
%
\special{pn 8}%
\special{pa 2952 830}%
\special{pa 2952 870}%
\special{da 0.070}%
\special{sh 1}%
\special{pa 2952 870}%
\special{pa 2972 804}%
\special{pa 2952 818}%
\special{pa 2932 804}%
\special{pa 2952 870}%
\special{fp}%
\put(27.1000,-9.7600){\makebox(0,0)[rb]{{\small $2v$}}}%
\put(44.3800,-14.5600){\makebox(0,0)[lt]{{\small $v$}}}%
\put(47.2900,-12.4500){\makebox(0,0)[lb]{{\small ${\bf J}_{\gamma(0)}(v)$}}}%
%
\special{pn 8}%
\special{pa 4690 1204}%
\special{pa 4322 1342}%
\special{da 0.070}%
\special{sh 1}%
\special{pa 4322 1342}%
\special{pa 4390 1338}%
\special{pa 4372 1324}%
\special{pa 4376 1300}%
\special{pa 4322 1342}%
\special{fp}%
\put(40.4900,-2.6800){\makebox(0,0)[lb]{{\small ${\rm Exp}_{\gamma(0)}({\bf J}_{\gamma(0)}(2v))$}}}%
\put(44.5700,-4.8800){\makebox(0,0)[lb]{{\small ${\rm Exp}_{\gamma(0)}({\bf J}_{\gamma(0)}(v))$}}}%
%
\special{pn 8}%
\special{ar 3312 1684 116 350  3.7714994 6.0839847}%
%
\special{pn 8}%
\special{pa 3224 1472}%
\special{pa 3224 1498}%
\special{fp}%
\special{sh 1}%
\special{pa 3224 1498}%
\special{pa 3244 1430}%
\special{pa 3224 1444}%
\special{pa 3204 1430}%
\special{pa 3224 1498}%
\special{fp}%
\put(32.1500,-16.5900){\makebox(0,0)[lt]{{\small $0$-section}}}%
%
\special{pn 8}%
\special{ar 4856 1652 118 350  3.7798148 6.0804038}%
%
\special{pn 8}%
\special{pa 4768 1448}%
\special{pa 4768 1464}%
\special{fp}%
\special{sh 1}%
\special{pa 4768 1464}%
\special{pa 4788 1398}%
\special{pa 4768 1412}%
\special{pa 4748 1398}%
\special{pa 4768 1464}%
\special{fp}%
\put(47.1900,-16.3500){\makebox(0,0)[lt]{{\small $G\cdot o(=G/K)$}}}%
%
\special{pn 8}%
\special{pa 4962 1562}%
\special{pa 4972 1620}%
\special{da 0.070}%
%
\special{pn 8}%
\special{pa 4876 748}%
\special{pa 4506 790}%
\special{da 0.070}%
\special{sh 1}%
\special{pa 4506 790}%
\special{pa 4574 802}%
\special{pa 4560 784}%
\special{pa 4570 762}%
\special{pa 4506 790}%
\special{fp}%
\put(49.1300,-7.9700){\makebox(0,0)[lb]{$\gamma_v^{\mathbb C}$}}%
%
\special{pn 8}%
\special{pa 4458 302}%
\special{pa 4312 976}%
\special{da 0.070}%
\special{sh 1}%
\special{pa 4312 976}%
\special{pa 4346 916}%
\special{pa 4324 924}%
\special{pa 4306 908}%
\special{pa 4312 976}%
\special{fp}%
%
\special{pn 8}%
\special{pa 4856 512}%
\special{pa 4312 1220}%
\special{da 0.070}%
\special{sh 1}%
\special{pa 4312 1220}%
\special{pa 4368 1180}%
\special{pa 4344 1178}%
\special{pa 4338 1156}%
\special{pa 4312 1220}%
\special{fp}%
%
\special{pn 8}%
\special{pa 2506 2034}%
\special{pa 3020 1758}%
\special{fp}%
%
\special{pn 8}%
\special{pa 2506 2222}%
\special{pa 3020 1920}%
\special{fp}%
%
\special{pn 8}%
\special{pa 4060 2230}%
\special{pa 4574 1936}%
\special{fp}%
\end{picture}%
\hspace{5.5truecm}}

\vspace{0.75truecm}

\centerline{{\bf Figure 2$\,$:$\,$ The image of $\gamma_{\ast}$}}

\vspace{0.35truecm}

\noindent
{\it Remark 2.1.} In \cite{K4}, this fact was stated but the proof was not given.  

\vspace{0.25truecm}

\noindent
{\bf 2.2.\ A construction of $G$-invariant K$\ddot{\rm a}$hler metrics on $G^{\mathbb C}/K^{\mathbb C}$}\ \ 
In this subsection, we recall a construction of $G$-invariant K$\ddot{\rm a}$hler metric on a complex symmetric space $G^{\mathbb C}/K^{\mathbb C}$.  
Let $G/K$ be a symmetric space of compact type and $(G^{\mathbb C}/K^{\mathbb C},\,{\bf J})$ be the complexification 
of $G/K$.  Let $\psi:G^{\mathbb C}/K^{\mathbb C}\to\mathbb R$ be a strictly plurisubharmonic function over 
$G^{\mathbb C}/K^{\mathbb C}$, where we note that ``strictly plurisubharmonicity'' means that 
the Hermitian matrix $\displaystyle{\left(\frac{\partial^2\psi}{\partial z_i\partial\bar z_j}\right)}$ 
is positive (or equivalently, $(\partial\overline{\partial}\psi)(\vZ,\overline \vZ)>0$ holds for any 
nonzero $(1,0)$-vector $\vZ$).  Then $\omega_{\psi}:=\sqrt{-1}\partial\overline{\partial}\psi$ 
is a real non-degenerate closed $2$-form on $G^{\mathbb C}/K^{\mathbb C}$ and the symmetric $(0,2)$-tensor field 
$\beta_{\psi}$ associated to ${\bf J}$ and $\omega_{\psi}$ is positive definite.  
Hence $({\bf J},\beta_{\psi})$ gives a K$\ddot{\rm a}$hler structure of $G^{\mathbb C}/K^{\mathbb C}$.  
Thus, from each strictly plurisubharmonic function over $G^{\mathbb C}/K^{\mathbb C}$, 
we can construct a K$\ddot{\rm a}$hler structure of $G^{\mathbb C}/K^{\mathbb C}$.  
Denote by ${\rm Exp}_p$ the exponential map of the anti-K$\ddot{\rm a}$hler manifold $(G^{\Bbb C}/K^{\Bbb C},\beta_A)$ 
at $p(\in G^{\Bbb C}/K^{\Bbb C})$ and $\exp$ the exponentional map of the Lie group $G^{\Bbb C}$.  
Set $\mathfrak g^d:=\mathfrak k\oplus\sqrt{-1}\mathfrak p(\subset\mathfrak g^{\Bbb C})$ and 
$G^d=\exp(\mathfrak g^d)$.  Denote by $\beta_{G/K}$ the $G$-invariant (Riemannian) metric on $G/K$ 
induced from $-B\vert_{\mathfrak p\times\mathfrak p}$ and $\beta_{G^d/K}$ the $G^d$-invariant negative definite metric 
on $G^d/K$ induced from $-({\rm Re}\,B^{\Bbb C})\vert_{\sqrt{-1}\mathfrak p\times\sqrt{-1}\mathfrak p}$.  
We may assume that the metric of $G/K$ is equal to $\beta_{G/K}$ by homothetically transforming the metric of 
$G/K$ if necessary.  On the other hand, the Riemannian manifold $(G^d/K,-\beta_{G^d/K})$ is a (Riemannian) 
symmetric space of non-compact type.  
The orbit $G\cdot o$ is isometric to $(G/K,\beta_{G/K})$ and the normal umbrella 
${\rm Exp}_o(T^{\perp}_o(G\cdot o))(=G^d\cdot o)$ is isometric to $(G^d/K,\beta_{G^d/K})$.  
The complexification $\mathfrak p^{\Bbb C}$ of $\mathfrak p$ is identified with $T_o(G^{\Bbb C}/K^{\Bbb C})$ and 
$\sqrt{-1}\mathfrak p$ is identified with $T_o({\rm Exp}_o(T^{\perp}_o(G\cdot o)))$.  
Let $\mathfrak a$ be a maximal abelian subspace of $\mathfrak p$, where we note that 
${\rm dim}\,\mathfrak a=r$.  Denote by $W$ the Weyl group of $G^d/K$ with 
respect to $\sqrt{-1}\mathfrak a$.  This group acts on $\sqrt{-1}\mathfrak a$.  
Let $C(\subset\sqrt{-1}\mathfrak a)$ be a Weyl domain (i.e., a fundamental domain of the action 
$W\curvearrowright\sqrt{-1}\mathfrak a$).  Then we have 
$G\cdot{\rm Exp}_o(\overline C)=G^{\Bbb C}/K^{\Bbb C}$, where $\overline C$ is the closure of $C$.  
For a $W$-invariant connected open neighborhood $D$ of $0$ in $\sqrt{-1}\mathfrak a$, we define a neighborhood 
$U_1(D)$ of $o$ 
in $\Sigma:={\rm Exp}_o(\sqrt{-1}\mathfrak a)$ by $U_1(D):={\rm Exp}_o(D)$, a neighborhood $U_2(D)$ of $o$ in 
$G^d/K$ by $U_2(D):=K\cdot U_1(D)$ and a tubular neighborhood $U_3(D)$ of $G\cdot o$ in $G^{\Bbb C}/K^{\Bbb C}$ by 
$U_3(D):=G\cdot U_1(D)$.  
Denote by ${\rm Conv}_W^+(D)$ the space of all $W$-invariant strictly convex ($C^{\infty}$-)functions over $D$, 
${\rm Conv}_K^+(U_2(D))$ the space of all $K$-invariant strictly convex ($C^{\infty}$-)functions 
over $U_2(D)$ and $PH_G^+(U_3(D))$ the space of all $G$-invariant strictly plurisubharmonic ($C^{\infty}$-)functions 
over $U_3(D)$, where ``strictly convex function'' means that its real Hessian is positive definite.  
The restriction map (which is  denoted by $\mathcal R_{32}^D$) from $U_3(D)$ to $U_2(D)$ gives 
an isomorphism of $PH_G^+(U_3(D))$ onto ${\rm Conv}_K^+(U_2(D))$ and the composition of the restriction map 
(which is denoted by $\mathcal R_{31}^D$) from $U_3(D)$ to $U_1(D)$ with ${\rm Exp}_o$ gives an isomorphism of 
$PH_G^+(U_3(D))$ onto ${\rm Conv}_W^+(D)$ (see \cite{AL}).  
Hence we suffice to construct $W$-invariant strictly convex functions over $D$ or $K$-invariant strictly convex 
functions over $U_2(D)$ to construct strictly plurisubharmonic functions over $U_3(D)$.  Let $\psi$ be a 
$G$-invariant strictly plurisubharmonic ($C^{\infty}$-)function over $U_3(D)$.  
Set $\bar{\psi}:=\mathcal R_{32}^D(\psi)$ and $\bar{\bar{\psi}}:=\mathcal R_{31}^D(\psi)\circ{\rm Exp}_o$.  
Conversely, for a $W$-invariant strictly convex ($C^{\infty}$-)function $\rho$ over $D$, 
denote by $\rho^h$ the $G$-invariant strictly plurisubharmonic ($C^{\infty}$-)function $\psi$ over $U_3(D)$ with 
$\bar{\bar{\psi}}=\rho$.  
Similarly, for a $K$-invariant strictly 
%
convex ($C^{\infty}$-)function $\sigma$ over $U_2(D)$, 
set $\overline{\sigma}:=\mathcal R_{21}^D(\sigma)\circ{\rm Exp}_o$ and denote by $\sigma^h$ the $G$-invariant 
strictly plurisubharmonic ($C^{\infty}$-)function $\psi$ over $U_3(D)$ with $\bar{\psi}=\sigma$.  
Also, for a $W$-invariant strictly convex ($C^{\infty}$-)function $\rho$ over $D$, 
denote by $\rho^d$ the $K$-invariant strictly convex $C^{\infty}$-function $\sigma$ over $U_2(D)$ with 
$\overline{\sigma}=\rho$.  
Denote by $Ric_{\psi}$ the Ricci form of $\beta_{\psi}$.  
It is known that $Ric_{\psi}$ is described as 
$$
Ric_{\psi}=-\sqrt{-1}\partial\overline{\partial}\log\,{\rm det}
\left(\frac{\partial^2\psi}{\partial z_i\partial\bar z_j}\right)\leqno{(2.1)}$$
(see P158 of \cite{KN}), where $(z_1,\cdots,z_n)$ is any complex coordinates of $G^{\Bbb C}/K^{\Bbb C}$.  
Note that the right-hand side of $(2.1)$ is independent of the choice of the complex coordinates 
$(z_1,\cdots,z_n)$ of $G^{\Bbb C}/K^{\Bbb C}$.  

\section{Proof of Theorem A} 
In this section, we shall prove Theorem A stated in Introduction.  
First we define a natural $G^{\mathbb C}$-invariant non-vanishing holomorphic $(n,0)$-form on $G^{\mathbb C}/K^{\mathbb C}$.  
Take an orthonormal basis $(\ve_1,\cdots,\ve_n)$ of $\mathfrak p$ with respect to $-B$ and let 
$(\theta^1,\cdots,\theta^n)$ be the dual basis of $(\ve_1,\cdots,\ve_n)$.  
Also, let $(\theta^i)^{\Bbb C}$ ($i=1,\cdots,n$) be the complexification of $\theta^i$.  
Since $(\theta^1)^{\Bbb C}\wedge\cdots\wedge(\theta^n)^{\Bbb C}$ is 
${\rm Ad}_{G^{\Bbb C}}(K^{\Bbb C})\vert_{\mathfrak p^{\Bbb C}}$-invariant, 
we obtain the $G^{\Bbb C}$-invariant non-vanishing holomorphic $(n,0)$-form $\Omega$ on $G^{\Bbb C}/K^{\Bbb C}$ with 
$\Omega_{eK^{\Bbb C}}=(\theta^1)^{\Bbb C}\wedge\cdots\wedge(\theta^n)^{\Bbb C}$.  

Let ${\bf J}$ be the $G^{\mathbb C}$-invariant complex structure with ${\bf J}_{eK}={\bf j}(=\sqrt{-1}\,{\rm id}_{\mathfrak p}$).  
Take any point $p_0:={\rm Exp}_o(\vZ_0)$ 
($\vZ_0\in\sqrt{-1}\mathfrak p$) 
and an orthonormal basis $(\ve_1,\cdots,\ve_n)$ of $\sqrt{-1}\mathfrak p$ with respect to 
$\beta_{G^d/K}$ satisfying $\ve_i\in\sqrt{-1}\mathfrak a$ ($i=1,\cdots,r$), where ${\rm Exp}_o$ is the exponential map of $(G^{\mathbb C}/K^{\mathbb C},\beta_A)$ at $o$.  
Let $\widetilde U_o$ be a sufficiently small neighborhood of the origin $0$ 
in $\mathfrak p^{\mathbb C}$.  
Define a local complex coordinates $\varphi_o=(z^o_1,\cdots,z^o_n)$ on 
$U_o:={\rm Exp}_o(\widetilde U_o)$ by 
$$({\rm Exp}_o\vert_{\widetilde U_o})^{-1}(p)=\sum_{i=1}^nz_i^o(p)\ve_i\quad\,\,(p\in U_o).$$
By Proposition 3.1 of \cite{K3}, $\varphi_o$ gives a local complex coordinates of the complex manifold $(G^{\mathbb C}/K^{\mathbb C},{\bf J})$.  
Set $U_{p_0}:=\exp(\vZ_0)(U_o)$, where $\exp$ is the exponential map of the Lie group $G^{\mathbb C}$.  
Define a local complex coordinates $\varphi^{p_0}=(z^{p_0}_1,\cdots,z^{p_0}_n)$ on $U_{p_0}$ by 
$z^{p_0}_i=z^o_i\circ\exp(\vZ_0)^{-1}$ ($i=1,\cdots,n$).  
Let $z^o_i=x^o_i+\sqrt{-1}y^o_i$ and $z^{p_0}_i=x^{p_0}_i+\sqrt{-1}y^{p_0}_i$ ($i=1,\cdots,n$).  
We call such a local complex coordinates $\varphi^{p_0}=(z^{p_0}_1,\cdots,z^{p_0}_n)$ 
{\it a normal complex coordinates about a point $p_0$ of the real form} $G^d\cdot o(=G^d/K)$ 
{\it associated to $(\ve_1,\cdots,\ve_n)$}.  
Set $\ve^{p_0}_i:=\exp_{G^d}(\vZ_0)_{\ast o}(\ve_i)$.  
Then we note that the following relation holds:
$${\rm Exp}_{p_0}^{-1}(p)=\sum_{i=1}^n\left(x^{p_0}_i(p)\ve^{p_0}_i+y^{p_0}_i(p)J_{p_0}(\ve^{p_0}_i)\right)
\quad(p\in U_{p_0}).$$

For the simplicity, we denote $\sqrt{-1}\mathfrak p$ and $\sqrt{-1}\mathfrak a$ by
$\mathfrak p^d$ and $\mathfrak a^d$, respectively.  
For $\lambda\in(\mathfrak a^d)^{\ast}$, we define $\mathfrak p^d_{\lambda}$ and $\mathfrak k_{\lambda}$ by 
$$\mathfrak k_{\lambda}:=\{v\in\mathfrak k\,\vert\,{\rm ad}(\vZ)^2(v)=\lambda(\vZ)^2v\,\,(\forall\,\vZ\in\mathfrak a^d)\}$$
and
$$\mathfrak p^d_{\lambda}:=\{v\in\mathfrak p^d\,\vert\,{\rm ad}(\vZ)^2(v)=\lambda(\vZ)^2v\,\,(\forall\,\vZ\in\mathfrak a^d)\}.$$
Also, we define $\triangle(\subset(\mathfrak a^d)^{\ast})$ by 
$$\triangle:=\{\lambda\in(\mathfrak a^d)^{\ast}\setminus\{0\}\,\vert\,\mathfrak p^d_{\lambda}\not=\{0\}\,\,\}.$$
Here we note that $\triangle$ is described as $\triangle=\triangle'\amalg\{2\lambda\,\vert\,\lambda\in S\}$ for some root system $\triangle'$ and some subset $S$ of $\triangle'$ 
(see \cite{H}).  However, we call $\triangle$ {\it the restricted root system associated to the symmetric pair} $(G^d,K)$ though $\triangle$ is not necessarily a root system.  
Let $\triangle_+$ be the positive root subsystem of $\triangle$ with respect to some lexicographic ordering of $(\mathfrak a^d)^{\ast}$.  
Then we have the following root space decompositions 
$$\mathfrak k=\mathfrak z_{\mathfrak k}(\mathfrak a^d)\oplus\left(
\mathop{\oplus}_{\lambda\in\triangle_+}\mathfrak k_{\lambda}\right)\quad{\rm and}\quad
\mathfrak p^d=\mathfrak a^d\oplus
\left(\mathop{\oplus}_{\lambda\in\triangle_+}\mathfrak p^d_{\lambda}\right),$$
where $\mathfrak z_{\mathfrak k}(\mathfrak a^d)$ is the centralizer of $\mathfrak a^d$ in $\mathfrak k$.  
Set $m_{\lambda}:={\rm dim}\,\mathfrak p^d_{\lambda}$ ($\lambda\in\triangle_+$).  
Note that the first root space decomposition is the simultaneously eigenspace decomposition of the commutative family 
$$\{{\rm ad}(\vZ)^2:\mathfrak k\to\mathfrak k\,\vert\,\vZ\in\mathfrak a^d\}$$
of the symmetric transformations of $\mathfrak k$ and that the second root space decomposition is the simultaneously eigenspace decomposition of the commutative family 
$$\{{\rm ad}(\vZ)^2:\mathfrak p^d\to\mathfrak p^d\,\vert\,\vZ\in\mathfrak a^d\}$$
of symmetric transformations of $\mathfrak p^d$, where ${\rm ad}$ is the adjoint representation of $\mathfrak g$.  

Let $\mathcal C(\subset\mathfrak a^d)$ be a Weyl domain (i.e., a fundamental domain of the Weyl group action 
$W\curvearrowright\mathfrak a^d$).  The Weyl domain $\mathcal C$ is given by 
$$\mathcal C:=\{\vv\in\mathfrak a^d\,\vert\,\lambda(\vv)>0\,\,(\forall\,\lambda\in\triangle_+)\}.$$
Points of $W\cdot\mathcal C$ and $G\cdot{\rm Exp}_o(W\cdot\mathcal C)$ are called {\it regular points}.  
We call $G\cdot{\rm Exp}_o(W\cdot\mathcal C)$ the {\it regular point set} of 
$G^{\mathbb C}/K^{\mathbb C}$ and denote it by $(G^{\mathbb C}/K^{\mathbb C})_{reg}$.  
Also, we call $K\cdot{\rm Exp}_o(W\cdot\mathcal C)$ the {\it regular point set} of 
$G^d/K$ and denote it by $(G^d/K)_{reg}$.  
Note that $(G^{\mathbb C}/K^{\mathbb C})_{reg}$ (resp. $(G^d/K)_{reg}$) is an open dense subset of 
$G^{\mathbb C}/K^{\mathbb C}$ (resp. $G^d/K$).  

Take a regular point $p:={\rm Exp}_o(\vZ)$ ($\vZ\in W\cdot\mathcal C$) of the isotropy group action $K\curvearrowright G^d/K$.  
Then we have 
$$(\exp\,\vZ)_{\ast}^{-1}(T_{p}(K\cdot p))=\mathop{\oplus}_{\lambda\in\triangle_+}\mathfrak p^d_{\lambda}\leqno{(3.1)}$$
and 
$$(\exp\,\vZ)_{\ast}^{-1}(T_{p}^{\perp}(K\cdot p))=\mathfrak a^d,$$
where $T_{p}^{\perp}(K\cdot p)$ is the normal space of $K\cdot p$ in $G^d\cdot o(=G^d/K)$.  
Denote by $A$ and $h$ the shape tensor and the second fundamental form of the principal orbit $K\cdot p$ (which is a submanifold in $G^d/K$), respectively. 
For any $\vv\in\mathfrak a^d$, we have 
$$(A_p)_{(\exp\,\vZ)_{\ast}(\vv)}\vert_{(\exp\,\vZ)_{\ast}(\mathfrak p^d_{\lambda})}=-\frac{\lambda(\vv)}{\tanh\lambda(\vZ)}\,{\rm id}\leqno{(3.2)}$$
(see \cite{V}).  

Both Lie groups $G$ and $K^{\mathbb C}$ are symmetric subgroups of $G^{\mathbb C}$, that is, 
the actions $G\curvearrowright G^{\mathbb C}/K^{\mathbb C}$ and $K^{\mathbb C}\curvearrowright G^{\mathbb C}/K^{\mathbb C}$ are Hermann type action, where we note that 
the terminology of ``Hermann type action'' was used in \cite{K2}.  
Let $\theta_G$ and $\theta_{K^{\mathbb C}}$ be the involutions of $G^{\mathbb C}$ with $({\rm Fix}\,\theta_G)_0\subset G\subset{\rm Fix}\,\theta_G$ and 
$({\rm Fix}\,\theta_{K^{\mathbb C}})_0\subset K^{\mathbb C}\subset{\rm Fix}\,\theta_{K^{\mathbb C}}$, respectively.  
Denote by the same symbols the differentials of $\theta_G$ and $\theta_{K^{\mathbb C}}$ at $e$, respectively.  
Note that $\theta_G$ is the Cartan involution of $G^{\mathbb C}$.  The eigenspace decompositions of $\theta_G$ and $\theta_{K^{\mathbb C}}$ are as 
$\mathfrak g^{\mathbb C}=\mathfrak g\oplus\sqrt{-1}\mathfrak g$ and $\mathfrak g^{\mathbb C}=\mathfrak k^{\mathbb C}\oplus\mathfrak p^{\mathbb C}$, respectively.  
The root space decomposition of $\mathfrak p^{\mathbb C}(=T_o(G^{\mathbb C}/K^{\mathbb C})$ for a maximal abelian subspace $\mathfrak a^d$ of 
$\mathfrak p^{\mathbb C}\cap\sqrt{-1}\mathfrak g=\mathfrak p^d$ is given by 
$$\mathfrak p^{\mathbb C}=\mathfrak a^{\mathbb C}\oplus\left(\mathop{\oplus}_{\lambda\in\triangle_+}\mathfrak p_{\lambda}^{\mathbb C}\right).$$
Note that this root space decomposition is the simultaneously eigenspace decomposition of the commutative family 
$$\{{\rm ad}(\vv)^2:\mathfrak p^{\mathbb C}\to\mathfrak p^{\mathbb C}\,\vert\,\vv\in\mathfrak a^d\}$$
of symmetric transformations of $\mathfrak p^{\mathbb C}$, where ${\rm ad}$ is the adjoint representation of $\mathfrak g^{\mathbb C}$.  
Also, the simultaneously eigenspace decomposition of the extended commutative family 
$$\{{\rm ad}(\vv)^2:\mathfrak p^{\mathbb C}\to\mathfrak p^{\mathbb C}\,\vert\,\vv\in\mathfrak a^d\}\cup\{\theta_G\vert_{\mathfrak p^{\mathbb C}}\}$$
of symmetric transformations of $\mathfrak p^{\mathbb C}$ is given by 
$$\mathfrak p^{\mathbb C}=\mathfrak a\oplus\mathfrak a^d\oplus\left(\mathop{\oplus}_{\lambda\in\triangle_+}\mathfrak p_{\lambda}\right)
\oplus\left(\mathop{\oplus}_{\lambda\in\triangle_+}\mathfrak p^d_{\lambda}\right).\leqno{(3.3)}$$

Fix a regular point $p:={\rm Exp}_o(\vZ)$ ($\vZ\in W\cdot\mathcal C$) of the action $G\curvearrowright G^{\mathbb C}/K^{\mathbb C}$.  
Then we have 
$$(\exp\,\vZ)_{\ast}^{-1}(T_{p}(G\cdot p))=\mathfrak a\oplus\left(\mathop{\oplus}_{\lambda\in\triangle_+}\mathfrak p_{\lambda}\right)\oplus
\left(\mathop{\oplus}_{\lambda\in\triangle_+}\mathfrak p^d_{\lambda}\right)\leqno{(3.4)}$$
and 
$$(\exp\,\vZ)_{\ast}^{-1}(T_{p}^{\perp}(G\cdot p))=\mathfrak a^d.$$
Denote by $\widehat A$ and $\widehat h$ the shape tensor and the second fundamental form of the regular orbit $G\cdot p$ (which is a submanifold in 
$G^{\mathbb C}/K^{\mathbb C}$), respectively.  
Then, according to $(6.1)$ and $(6.2)$ in \cite{K2}, we have 
$$(\widehat A_p)_{(\exp\,\vZ)_{\ast}(\vv)}\vert_{(\exp\,\vZ)_{\ast}(\mathfrak p^d_{\lambda})}=-\frac{\lambda(\vv)}{\tanh\lambda(\vZ)}\,{\rm id}\quad\,\,
(\lambda\in\triangle_+)\leqno{(3.5)}$$
and 
$$(\widehat A_p)_{(\exp\,\vZ)_{\ast}(\vv)}\vert_{(\exp\,\vZ)_{\ast}(\mathfrak p_{\lambda})}=-\lambda(\vv)\tanh\lambda(\vZ)\,{\rm id}\quad\,\,(\lambda\in\triangle_+)
\leqno{(3.6)}$$
for $\vv\in(\exp\,\vZ)_{\ast}^{-1}(T_{p}^{\perp}(G\cdot o))$.  

Take an orthonormal base $(\ve^0_i)_{i=1}^r$ of $\mathfrak a^d$ and an orthonormal base $(\ve^{\lambda}_i)_{i=1}^{m_{\lambda}}$ of $\mathfrak p^d_{\lambda}$ with respect to 
$\beta_A$.  Let $\triangle_+=\{\lambda_1,\cdots,\lambda_l\}$ ($\lambda_1<\lambda_2<\cdots<\lambda_l$), where $<$ is a lexicographical ordering 
of $\mathfrak a^d$ with respect to a basis of $\mathfrak a^d$.  
For the simplicity, set $m_0:=r$ and $m_i:=r+m_{\lambda_1}+\cdots+m_{\lambda_i}$ ($i=1,\cdots,l$).  
Define $(\ve_1,\cdots,\ve_n)$ by 
$$\ve_i:=\left\{\begin{array}{ll}
\ve^0_i & (1\leq i\leq m_0)\\
\ve^{\lambda_1}_{i-m_0} & (m_0+1\leq i\leq m_1)\\
\ve^{\lambda_2}_{i-m_1} & (m_1+1\leq i\leq m_2)\\
\vdots & \qquad\qquad\qquad\qquad\,\,\vdots \\
\ve^{\lambda_l}_{i-m_{l-1}} & (m_{l-1}+1\leq i\leq m_l)
\end{array}\right.\leqno{(3.7)}$$
and $(\ve^{p}_1,\cdots,\ve^{p}_n)$ by $\ve^{p}_i:=(\exp\,\vZ)_{\ast}(\ve_i)$ ($i=1,\cdots,n$).  
For the convenience, set $\mathfrak p^d_0:=\mathfrak a^d$ and $\lambda_0:=0$.  

\vspace{0.5truecm}

\noindent
{\bf Lemma 3.1.} {\sl Assume that $p\in(G^d/K)_{\rm reg}$ (i.e., $\vZ\in\mathcal C$).  Let $(U,\varphi=(z_1,\cdots,z_n))$ be a normal complex coordinates about 
a point $p$ of the real form $G^d\cdot o(=G^d/K)$ associated to $(\ve^{p}_1,\cdots,\ve^{p}_n)$.  
Then, for any $\rho\in{\rm Conv}_W^+(D)$, we have 
{\small
$$\left.\frac{\partial^2\rho^h}{\partial z_i\partial\overline z_j}\right\vert_{p}
=\left\{\begin{array}{ll}
\displaystyle{\frac{1}{4}\,(\nabla^dd\rho^d)_{p}((\ve_i^0)^p,(\ve_i^0)^p)} & (i=j\in\{1,\cdots,m_0\})\\
\displaystyle{\frac{1}{4}\,(\,\,(\nabla^dd\rho^d)_{p}((\ve_{i-m_0}^{\lambda_1})^p,(\ve_{i-m_0}^{\lambda_1})^p)} & (i=j\in\{m_0+1,\cdots\\
\displaystyle{\qquad +\tanh\lambda_1(\vZ)\cdot\lambda_1(({\rm grad}\,\rho)_{\vZ})\delta_{ij}\,\,)} & \hspace{1.15truecm} \cdots,m_1\})\\
\qquad\qquad\qquad\qquad\,\,\vdots & \qquad\qquad\qquad\qquad\,\,\vdots\\
\displaystyle{\frac{1}{4}\,(\,\,(\nabla^dd\rho^d)_{p}((\ve^{\lambda_l}_{i-m_{l-1}})^p,(\ve^{\lambda_l}_{i-m_{l-1}})^p)} & (i=j\in\{m_{l-1}+1,\cdots\\
\hspace{0.55truecm}\displaystyle{+\tanh\,\lambda_l(\vZ)\cdot\lambda_l(({\rm grad}\,\rho)_{\vZ})\delta_{ij})} & \hspace{1.15truecm}\cdots,m_l\})\\
\displaystyle{\frac{1}{4}\,(\nabla^dd\rho^d)_{p}(\ve_i^{p},\ve_j^{p})} & (i\not=j).
\end{array}\right.
$$
}
where $\lambda_a^{\sharp}$ denotes the element of $(\exp\,\vZ)_{\ast}^{-1}(T_{p}^{\perp}(G\cdot p))$ such that 
$\beta_A(\lambda_a^{\sharp},\vv)=\lambda_a(\vv)$ holds for any $\vv\in(\exp\,\vZ)_{\ast}^{-1}(T_{p}^{\perp}(G\cdot o))$.  
}

\vspace{0.5truecm}

\noindent
{\it Proof.}\ \ 
First we note that {\small$({\bf J}\ve^p_1,\cdots,{\bf J}\ve^p_n)$} be an orthonormal basis of the normal space 
$T_{p}^{\perp}(G^d/K)\,(\subset T_{p}(G\cdot p))$.  

Let $i,j\in\{1,\cdots,m_0\}$.  
Define a map $\delta^0:[0,\infty)^2\to{\rm Exp}_p((\exp\,Z)_{\ast}(\mathfrak a^{\mathbb C}))$ by 
$$\delta^0(s,t):={\rm Exp}_p((\exp\,Z)_{\ast}(t\ve^p_i+\sqrt{-1}s\ve^p_j)\quad\,\,((s,t)\in\mathbb R^2)$$
(see Figure 3).  
Then, since $\rho^h$ is $G$-invariant and ${\rm Exp}_p((\exp\,Z)_{\ast}(\mathfrak a^{\mathbb C}))$ is totally geodesic in $G^{\mathbb C}/K^{\mathbb  C}$, 
we can show that 
$$(\nabla d\rho^h)_{p}(\ve^{p}_i,{\bf J}\ve^{p}_j)
=\left.\frac{\partial^2(\rho^h\circ\delta^0)}{\partial t\partial s}\right\vert_{s=0}-(d\rho^h)_{p}(\nabla_{\ve^p_i}
\frac{\partial\delta^0}{\partial s})={\bf 0}.$$
Similarly, we have 
$$(\nabla d\rho^h)_{p}({\bf J}\ve^{p}_i,{\bf J}\ve^{p}_j)={\bf 0}.$$
Therefore, we obtain 
\begin{align*}
\left.\frac{\partial^2\rho^h}{\partial z_i\partial\overline z_j}\right\vert_{p}
&=(\nabla d\rho^h)_p^{\mathbb C}\left(\left(\frac{\partial}{\partial z_i}\right)_p,\,\left(\frac{\partial}{\partial\overline z_j}\right)_p\right)\\
&=\frac{1}{4}(\nabla d\rho^h)_p^{\mathbb C}\left((\ve^0_i)^p-\sqrt{-1}{\bf J}(\ve^0_i)^p,\,(\ve^0_j)^p+\sqrt{-1}{\bf J}(\ve^0_j)^p\right)\\
&=\frac{1}{4}\,(\nabla^dd\rho^d)_{p}((\ve_i^0)^p,(\ve_j^0)^p).
\end{align*}

Let $i,j\geq m_0+1$.  
Let 
$\widehat{\gamma}_i$ be the geodesic in $G\cdot p$ with $\widehat{\gamma}_i'(0)={\bf J}\ve^{p}_i$.  
Then, from the $G$-invariance of $\rho^h$, we can show that 
\begin{align*}
(\nabla d\rho^h)_{p}({\bf J}\ve^{p}_i,{\bf J}\ve^{p}_i)
&=\left.\frac{d^2(\rho^h\circ\widehat{\gamma}_i)}{ds^2}\right\vert_{s=0}-(d\rho^h)_{p}(\nabla_{\widehat{\gamma}_i'(0)}\widehat{\gamma}_i')\\
&=-(d\rho^d)_{p}(\widehat h_p({\bf J}\ve^{p}_i,{\bf J}\ve^{p}_i)).
\end{align*}
Similarly, we can show 
$$(\nabla d\rho^h)_{p}({\bf J}e^{p}_i+{\bf J}e^{p}_j,{\bf J}e^{p}_i+{\bf J}e^{p}_j)=
-(d\rho^d)_{p}(\widehat h_p({\bf J}\ve^{p}_i+{\bf J}\ve^{p}_j,{\bf J}\ve^{p}_i+{\bf J}\ve^{p}_j)).$$
Hence, from the symmetricnesses of $(\nabla d\rho^h)_{p}$ and $\widehat h_{p}$, we have 
$$(\nabla d\rho^h)_{p}({\bf J}e^{p}_i,{\bf J}e^{p}_j)=-(d\rho^d)_{p}(\widehat h_p({\bf J}\ve^{p}_i,{\bf J}\ve^{p}_j)).\leqno{(3.8)}$$
On the other hand, according to $(3.6)$, we have 
$$\widehat h_p({\bf J}\ve^p_i,\,{\bf J}e^p_j)=
\left\{\begin{array}{ll}
\displaystyle{-\delta_{ij}\,\tanh\lambda_a(\vZ)\,\lambda_a^{\sharp}} & 
\displaystyle{\left(\begin{array}{r}
(i,j)\in\{m_{a-1}+1,\cdots,m_a\}^2\\
(a=1,\cdots,l)
\end{array}\right)}\\
0 & ((i,j):\,{\rm other})
\end{array}\right.
$$
($\delta_{ij}\,:\,$ the Kronecker's delta), where we used also $\beta_A({\bf J}\ve^p_i,\,{\bf J}\ve^p_j)=\delta_{ij}$.
Hence we obtain 
{\small
$$\begin{array}{l}
\hspace{0.5truecm}\displaystyle{(\nabla d\rho^h)_{p}({\bf J}e^p_i,\,{\bf J}e^p_j)}\\
\displaystyle{=
\left\{\begin{array}{ll}
\delta_{ij}\,\tanh\lambda_a(\vZ)\,\lambda_a(({\rm grad}\,\rho)_{\vZ}) & 
\displaystyle{\left(\begin{array}{r}
(i,j)\in\{m_{a-1}+1,\cdots,m_a\}^2\\
(a=1,\cdots,l)
\end{array}\right)}\\
0 & ((i,j):\,{\rm other}).
\end{array}\right.}
\end{array}\leqno{(3.9)}$$
}

Let $i,j\geq m_0+1$.  
Let $\gamma_i$ be the geodesic in $K\cdot p$ with $\gamma_i'(0)=\ve^{p}_i$ and 
$\widetilde{{\bf J}\ve_j^p}$ be the parallel normal vector field of 
the submanifold $K\cdot p$ in $G\cdot p$ 
along $\gamma_i$ with $(\widetilde{{\bf J}\ve^p_j})_0={\bf J}\ve^p_j$.  
Also, let $\widehat{\gamma}_j^t$ be the geodesic in $G\cdot p$ 
with $(\widehat{\gamma}_j^t)'(0)=(\widetilde{{\bf J}\ve^p_j})_t$.  
Define a geodesic variation $\delta$ by $\delta(s,t):=\widehat{\gamma}_j^t(s)$ (see Figure 4).  
Then, from the $G$-invariance of $\rho^h$, we have 
$$\begin{array}{l}
\displaystyle{(\nabla d\rho^h)_p({\bf J}\ve^p_i,\,\ve^p_j)=\left.\frac{\partial^2(\rho^h\circ\delta)}{\partial s\partial t}\right\vert_{s=t=0}
-(d\rho^h)_p(\nabla_{\ve^p_j}\widetilde{{\bf J}\ve^p_i})}\\
\hspace{2.85truecm}
\displaystyle{=-(d\rho^d)_{p}(\widehat h_p(\ve^{p}_j,{\bf J}\ve^{p}_i))=0.}
\end{array}\leqno{(3.10)}$$

\vspace{0.5truecm}

\centerline{
\unitlength 0.1in
\begin{picture}( 42.1500, 32.0000)(  7.2300,-39.3000)
%
\special{pn 20}%
\special{sh 1}%
\special{ar 2478 2770 10 10 0  6.28318530717959E+0000}%
\special{sh 1}%
\special{ar 2478 2770 10 10 0  6.28318530717959E+0000}%
%
\special{pn 8}%
\special{pa 1990 2430}%
\special{pa 1862 2430}%
\special{fp}%
%
\special{pn 8}%
\special{pa 1972 3072}%
\special{pa 1842 3072}%
\special{fp}%
\put(32.8100,-9.0000){\makebox(0,0)[lb]{$G^d\cdot o$}}%
\put(33.8000,-11.4600){\makebox(0,0)[lb]{${\rm Exp}_o(\mathfrak a^d)$}}%
\put(44.2200,-25.1400){\makebox(0,0)[lb]{$G\cdot o$}}%
%
\special{pn 8}%
\special{pa 4502 976}%
\special{pa 4938 976}%
\special{fp}%
%
\special{pn 8}%
\special{pa 4512 976}%
\special{pa 4512 760}%
\special{fp}%
\put(45.6100,-9.4800){\makebox(0,0)[lb]{$G^{\mathbb C}/K^{\mathbb C}$}}%
%
\special{pn 13}%
\special{pa 1990 2760}%
\special{pa 1862 2760}%
\special{fp}%
%
\special{pn 20}%
\special{sh 1}%
\special{ar 2478 2430 10 10 0  6.28318530717959E+0000}%
\special{sh 1}%
\special{ar 2478 2430 10 10 0  6.28318530717959E+0000}%
\put(40.7500,-20.6100){\makebox(0,0)[lb]{{\small This inside tube is $G\cdot p$.}}}%
\put(31.1200,-33.6400){\makebox(0,0)[lt]{$o$}}%
%
\special{pn 8}%
\special{pa 3520 1562}%
\special{pa 3102 1760}%
\special{fp}%
\special{sh 1}%
\special{pa 3102 1760}%
\special{pa 3172 1748}%
\special{pa 3150 1736}%
\special{pa 3154 1712}%
\special{pa 3102 1760}%
\special{fp}%
%
\special{pn 8}%
\special{ar 3212 2750 130 316  0.1664668 6.1359336}%
%
\special{pn 13}%
\special{ar 2468 2750 128 324  0.1696249 6.1358424}%
%
\special{pn 8}%
\special{pa 2894 806}%
\special{pa 2060 1364}%
\special{pa 2060 3930}%
\special{pa 2894 3346}%
\special{pa 2894 3346}%
\special{pa 2894 806}%
\special{fp}%
%
\special{pn 8}%
\special{pa 2468 3656}%
\special{pa 2458 1108}%
\special{fp}%
%
\special{pn 20}%
\special{pa 2458 2430}%
\special{pa 2458 1608}%
\special{fp}%
%
\special{pn 13}%
\special{pa 2498 2430}%
\special{pa 3222 2430}%
\special{fp}%
%
\special{pn 13}%
\special{pa 2478 2242}%
\special{pa 3202 2242}%
\special{fp}%
%
\special{pn 13}%
\special{pa 2460 2130}%
\special{pa 3214 2130}%
\special{fp}%
%
\special{pn 13}%
\special{pa 2460 1870}%
\special{pa 3204 1870}%
\special{fp}%
%
\special{pn 13}%
\special{pa 2478 1638}%
\special{pa 3202 1638}%
\special{fp}%
%
\special{pn 8}%
\special{pa 3062 1448}%
\special{pa 2478 1760}%
\special{fp}%
\special{sh 1}%
\special{pa 2478 1760}%
\special{pa 2546 1746}%
\special{pa 2524 1734}%
\special{pa 2526 1710}%
\special{pa 2478 1760}%
\special{fp}%
\put(31.1200,-14.5700){\makebox(0,0)[lb]{$\gamma_i$}}%
\put(35.6900,-15.5200){\makebox(0,0)[lb]{$\delta^0$}}%
\put(16.2300,-24.6700){\makebox(0,0)[rb]{$K\cdot p$}}%
%
\special{pn 20}%
\special{sh 1}%
\special{ar 2458 2420 10 10 0  6.28318530717959E+0000}%
\special{sh 1}%
\special{ar 2458 2420 10 10 0  6.28318530717959E+0000}%
\put(19.1100,-21.1800){\makebox(0,0)[rb]{$p$}}%
%
\special{pn 8}%
\special{pa 1962 2108}%
\special{pa 2458 2420}%
\special{fp}%
\special{sh 1}%
\special{pa 2458 2420}%
\special{pa 2412 2368}%
\special{pa 2412 2392}%
\special{pa 2390 2402}%
\special{pa 2458 2420}%
\special{fp}%
%
\special{pn 8}%
\special{pa 1664 2440}%
\special{pa 2336 2646}%
\special{fp}%
\special{sh 1}%
\special{pa 2336 2646}%
\special{pa 2278 2608}%
\special{pa 2284 2630}%
\special{pa 2266 2646}%
\special{pa 2336 2646}%
\special{fp}%
%
\special{pn 8}%
\special{pa 3320 1080}%
\special{pa 2478 1468}%
\special{fp}%
\special{sh 1}%
\special{pa 2478 1468}%
\special{pa 2546 1458}%
\special{pa 2526 1446}%
\special{pa 2530 1422}%
\special{pa 2478 1468}%
\special{fp}%
%
\special{pn 8}%
\special{pa 3222 872}%
\special{pa 2736 1128}%
\special{fp}%
\special{sh 1}%
\special{pa 2736 1128}%
\special{pa 2804 1114}%
\special{pa 2782 1102}%
\special{pa 2786 1078}%
\special{pa 2736 1128}%
\special{fp}%
%
\special{pn 13}%
\special{pa 2100 2750}%
\special{pa 2230 2760}%
\special{fp}%
%
\special{pn 13}%
\special{pa 2398 2760}%
\special{pa 2964 2770}%
\special{fp}%
%
\special{pn 13}%
\special{pa 3122 2760}%
\special{pa 4224 2770}%
\special{fp}%
%
\special{pn 8}%
\special{pa 1862 2232}%
\special{pa 2020 2232}%
\special{fp}%
%
\special{pn 8}%
\special{pa 2110 2232}%
\special{pa 4184 2232}%
\special{fp}%
%
\special{pn 8}%
\special{pa 1842 3222}%
\special{pa 2000 3222}%
\special{fp}%
%
\special{pn 8}%
\special{pa 2120 3222}%
\special{pa 4194 3222}%
\special{fp}%
%
\special{pn 8}%
\special{pa 4354 2506}%
\special{pa 4046 2760}%
\special{fp}%
\special{sh 1}%
\special{pa 4046 2760}%
\special{pa 4110 2734}%
\special{pa 4086 2726}%
\special{pa 4084 2702}%
\special{pa 4046 2760}%
\special{fp}%
%
\special{pn 8}%
\special{pa 3092 3318}%
\special{pa 2498 2780}%
\special{fp}%
\special{sh 1}%
\special{pa 2498 2780}%
\special{pa 2534 2840}%
\special{pa 2538 2816}%
\special{pa 2560 2810}%
\special{pa 2498 2780}%
\special{fp}%
\put(24.1700,-21.9300){\makebox(0,0)[rb]{$\ve^p_i$}}%
%
\special{pn 20}%
\special{pa 2458 2420}%
\special{pa 2458 2128}%
\special{fp}%
\special{sh 1}%
\special{pa 2458 2128}%
\special{pa 2438 2194}%
\special{pa 2458 2180}%
\special{pa 2478 2194}%
\special{pa 2458 2128}%
\special{fp}%
%
\special{pn 8}%
\special{ar 4016 2572 278 586  3.1580747 4.7389080}%
%
\special{pn 8}%
\special{pa 3738 2552}%
\special{pa 3738 2590}%
\special{fp}%
\special{sh 1}%
\special{pa 3738 2590}%
\special{pa 3758 2524}%
\special{pa 3738 2538}%
\special{pa 3718 2524}%
\special{pa 3738 2590}%
\special{fp}%
%
\special{pn 8}%
\special{ar 3996 2354 378 614  3.2309896 4.7123890}%
%
\special{pn 8}%
\special{pa 3628 2288}%
\special{pa 3628 2316}%
\special{fp}%
\special{sh 1}%
\special{pa 3628 2316}%
\special{pa 3648 2250}%
\special{pa 3628 2264}%
\special{pa 3608 2250}%
\special{pa 3628 2316}%
\special{fp}%
\put(40.4500,-18.0600){\makebox(0,0)[lb]{{\small This outside tube is $G\cdot\gamma_i(3/4)$.}}}%
%
\special{pn 13}%
\special{pa 2770 2440}%
\special{pa 2770 1648}%
\special{fp}%
%
\special{pn 13}%
\special{pa 3080 2420}%
\special{pa 3080 1628}%
\special{fp}%
%
\special{pn 13}%
\special{pa 3212 2430}%
\special{pa 3212 1638}%
\special{fp}%
%
\special{pn 8}%
\special{ar 3202 2732 200 492  0.1448125 6.2705960}%
%
\special{pn 8}%
\special{ar 2458 2722 198 492  0.1448125 6.2687978}%
%
\special{pn 13}%
\special{pa 2278 2760}%
\special{pa 2318 2760}%
\special{fp}%
%
\special{pn 8}%
\special{pa 2090 3072}%
\special{pa 2288 3072}%
\special{fp}%
%
\special{pn 8}%
\special{pa 2358 3072}%
\special{pa 3034 3072}%
\special{fp}%
%
\special{pn 8}%
\special{pa 3102 3072}%
\special{pa 4184 3072}%
\special{fp}%
%
\special{pn 8}%
\special{pa 3082 2430}%
\special{pa 4164 2430}%
\special{fp}%
%
\special{pn 8}%
\special{pa 2090 2430}%
\special{pa 2288 2430}%
\special{fp}%
%
\special{pn 8}%
\special{pa 2328 2430}%
\special{pa 3004 2430}%
\special{fp}%
%
\special{pn 13}%
\special{pa 3034 2770}%
\special{pa 3072 2770}%
\special{fp}%
%
\special{pn 20}%
\special{pa 2480 2430}%
\special{pa 2760 2440}%
\special{fp}%
\special{sh 1}%
\special{pa 2760 2440}%
\special{pa 2694 2418}%
\special{pa 2708 2438}%
\special{pa 2694 2458}%
\special{pa 2760 2440}%
\special{fp}%
%
\special{pn 8}%
\special{pa 3350 1910}%
\special{pa 2640 2410}%
\special{fp}%
\special{sh 1}%
\special{pa 2640 2410}%
\special{pa 2706 2388}%
\special{pa 2684 2380}%
\special{pa 2684 2356}%
\special{pa 2640 2410}%
\special{fp}%
\put(33.3000,-18.9000){\makebox(0,0)[lb]{${\bf J}\ve^p_j$}}%
\end{picture}%
\hspace{4truecm}}

\vspace{0.35truecm}

\centerline{{\bf Figure 3 $\,:\,$ Geodesic variation $\delta^0$}}

\vspace{0.75truecm}

\centerline{
\unitlength 0.1in
\begin{picture}( 42.7400, 29.2100)(  7.1600,-39.9100)
%
\special{pn 20}%
\special{sh 1}%
\special{ar 2482 2836 10 10 0  6.28318530717959E+0000}%
\special{sh 1}%
\special{ar 2482 2836 10 10 0  6.28318530717959E+0000}%
%
\special{pn 8}%
\special{pa 1994 2488}%
\special{pa 1866 2488}%
\special{fp}%
%
\special{pn 8}%
\special{pa 1974 3146}%
\special{pa 1846 3146}%
\special{fp}%
\put(32.7800,-13.3200){\makebox(0,0)[lb]{$G^d\cdot o$}}%
\put(33.7800,-15.8400){\makebox(0,0)[lb]{${\rm Exp}_o(\mathfrak a^d)$}}%
\put(44.3300,-25.7400){\makebox(0,0)[lb]{$G\cdot o$}}%
%
\special{pn 8}%
\special{pa 4552 1294}%
\special{pa 4990 1294}%
\special{fp}%
%
\special{pn 8}%
\special{pa 4562 1294}%
\special{pa 4562 1070}%
\special{fp}%
\put(46.1200,-12.6400){\makebox(0,0)[lb]{$G^{\mathbb C}/K^{\mathbb C}$}}%
%
\special{pn 13}%
\special{pa 2020 2840}%
\special{pa 1892 2840}%
\special{fp}%
%
\special{pn 20}%
\special{sh 1}%
\special{ar 2482 2488 10 10 0  6.28318530717959E+0000}%
\special{sh 1}%
\special{ar 2482 2488 10 10 0  6.28318530717959E+0000}%
\put(40.8400,-21.0800){\makebox(0,0)[lb]{{\small This tube is $G\cdot p$.}}}%
\put(31.1900,-34.4700){\makebox(0,0)[lt]{$o$}}%
%
\special{pn 8}%
\special{ar 3218 2818 130 326  0.1686528 6.1342954}%
%
\special{pn 13}%
\special{ar 2472 2818 128 332  0.1691665 6.1338206}%
%
\special{pn 8}%
\special{pa 2890 1196}%
\special{pa 2054 1694}%
\special{pa 2054 3992}%
\special{pa 2890 3468}%
\special{pa 2890 3468}%
\special{pa 2890 1196}%
\special{fp}%
%
\special{pn 13}%
\special{pa 2502 2488}%
\special{pa 3228 2488}%
\special{fp}%
%
\special{pn 20}%
\special{pa 2472 2478}%
\special{pa 2712 2478}%
\special{fp}%
\special{sh 1}%
\special{pa 2712 2478}%
\special{pa 2644 2458}%
\special{pa 2658 2478}%
\special{pa 2644 2498}%
\special{pa 2712 2478}%
\special{fp}%
\put(32.7800,-18.7500){\makebox(0,0)[lb]{$\delta$}}%
\put(16.2600,-31.2700){\makebox(0,0)[rt]{$K\cdot p$}}%
%
\special{pn 20}%
\special{sh 1}%
\special{ar 2462 2478 10 10 0  6.28318530717959E+0000}%
\special{sh 1}%
\special{ar 2462 2478 10 10 0  6.28318530717959E+0000}%
\put(19.0500,-15.9400){\makebox(0,0)[rb]{$p$}}%
%
\special{pn 8}%
\special{pa 1666 2496}%
\special{pa 2340 2710}%
\special{fp}%
\special{sh 1}%
\special{pa 2340 2710}%
\special{pa 2282 2670}%
\special{pa 2290 2694}%
\special{pa 2270 2708}%
\special{pa 2340 2710}%
\special{fp}%
%
\special{pn 20}%
\special{pa 2810 3972}%
\special{pa 3050 3972}%
\special{fp}%
\special{sh 1}%
\special{pa 3050 3972}%
\special{pa 2982 3952}%
\special{pa 2996 3972}%
\special{pa 2982 3992}%
\special{pa 3050 3972}%
\special{fp}%
\put(30.9900,-38.4500){\makebox(0,0)[lt]{'s $:\,\widetilde{{\bf J}\ve_j^p}$}}%
%
\special{pn 8}%
\special{pa 3318 1516}%
\special{pa 2472 1914}%
\special{fp}%
\special{sh 1}%
\special{pa 2472 1914}%
\special{pa 2542 1904}%
\special{pa 2520 1892}%
\special{pa 2524 1868}%
\special{pa 2472 1914}%
\special{fp}%
%
\special{pn 8}%
\special{pa 3218 1304}%
\special{pa 2732 1566}%
\special{fp}%
\special{sh 1}%
\special{pa 2732 1566}%
\special{pa 2800 1552}%
\special{pa 2778 1540}%
\special{pa 2780 1516}%
\special{pa 2732 1566}%
\special{fp}%
%
\special{pn 13}%
\special{pa 2080 2830}%
\special{pa 2296 2844}%
\special{fp}%
%
\special{pn 13}%
\special{pa 2402 2826}%
\special{pa 3070 2842}%
\special{fp}%
%
\special{pn 13}%
\special{pa 3130 2830}%
\special{pa 4236 2840}%
\special{fp}%
%
\special{pn 8}%
\special{pa 4364 2564}%
\special{pa 4056 2826}%
\special{fp}%
\special{sh 1}%
\special{pa 4056 2826}%
\special{pa 4120 2798}%
\special{pa 4096 2792}%
\special{pa 4094 2768}%
\special{pa 4056 2826}%
\special{fp}%
%
\special{pn 8}%
\special{pa 3100 3400}%
\special{pa 2502 2846}%
\special{fp}%
\special{sh 1}%
\special{pa 2502 2846}%
\special{pa 2538 2906}%
\special{pa 2542 2882}%
\special{pa 2564 2878}%
\special{pa 2502 2846}%
\special{fp}%
\put(18.2500,-20.4000){\makebox(0,0)[rb]{$\ve^p_i$}}%
%
\special{pn 8}%
\special{ar 4026 2632 278 602  3.1580145 4.7389080}%
%
\special{pn 8}%
\special{pa 3746 2614}%
\special{pa 3746 2652}%
\special{fp}%
\special{sh 1}%
\special{pa 3746 2652}%
\special{pa 3766 2586}%
\special{pa 3746 2600}%
\special{pa 3726 2586}%
\special{pa 3746 2652}%
\special{fp}%
%
\special{pn 8}%
\special{pa 2090 3150}%
\special{pa 4170 3150}%
\special{fp}%
%
\special{pn 8}%
\special{pa 3090 2488}%
\special{pa 4174 2488}%
\special{fp}%
%
\special{pn 8}%
\special{pa 2094 2488}%
\special{pa 2294 2488}%
\special{fp}%
%
\special{pn 8}%
\special{pa 2334 2488}%
\special{pa 3010 2488}%
\special{fp}%
%
\special{pn 8}%
\special{ar 2472 2818 10 534  0.0000000 6.2831853}%
%
\special{pn 20}%
\special{ar 2472 2826 130 330  3.3336407 4.7123890}%
%
\special{pn 8}%
\special{pa 1676 3146}%
\special{pa 2344 2952}%
\special{fp}%
\special{sh 1}%
\special{pa 2344 2952}%
\special{pa 2274 2952}%
\special{pa 2292 2968}%
\special{pa 2286 2990}%
\special{pa 2344 2952}%
\special{fp}%
\put(16.1600,-25.1600){\makebox(0,0)[rb]{$\gamma_i$}}%
%
\special{pn 13}%
\special{pa 2412 2574}%
\special{pa 3140 2574}%
\special{fp}%
%
\special{pn 13}%
\special{pa 2382 2682}%
\special{pa 3110 2682}%
\special{fp}%
%
\special{pn 13}%
\special{pa 2352 2768}%
\special{pa 3080 2768}%
\special{fp}%
%
\special{pn 20}%
\special{pa 2462 2488}%
\special{pa 2294 2584}%
\special{fp}%
\special{sh 1}%
\special{pa 2294 2584}%
\special{pa 2362 2568}%
\special{pa 2340 2558}%
\special{pa 2342 2534}%
\special{pa 2294 2584}%
\special{fp}%
%
\special{pn 8}%
\special{pa 1866 2070}%
\special{pa 2382 2526}%
\special{fp}%
\special{sh 1}%
\special{pa 2382 2526}%
\special{pa 2346 2466}%
\special{pa 2342 2490}%
\special{pa 2320 2496}%
\special{pa 2382 2526}%
\special{fp}%
%
\special{pn 8}%
\special{pa 1906 1614}%
\special{pa 2462 2478}%
\special{fp}%
\special{sh 1}%
\special{pa 2462 2478}%
\special{pa 2444 2410}%
\special{pa 2434 2432}%
\special{pa 2410 2432}%
\special{pa 2462 2478}%
\special{fp}%
%
\special{pn 13}%
\special{ar 2612 2818 130 330  3.3356988 4.7123890}%
%
\special{pn 13}%
\special{ar 2780 2818 130 330  3.3361040 4.7123890}%
%
\special{pn 13}%
\special{ar 2920 2826 130 330  3.3336407 4.7123890}%
%
\special{pn 13}%
\special{ar 3050 2818 130 330  3.3356988 4.7123890}%
%
\special{pn 13}%
\special{ar 3160 2818 130 330  3.3336407 4.7123890}%
%
\special{pn 8}%
\special{ar 3358 2546 558 728  3.1611469 4.5356801}%
%
\special{pn 8}%
\special{pa 2800 2526}%
\special{pa 2800 2546}%
\special{fp}%
\special{sh 1}%
\special{pa 2800 2546}%
\special{pa 2820 2478}%
\special{pa 2800 2492}%
\special{pa 2780 2478}%
\special{pa 2800 2546}%
\special{fp}%
%
\special{pn 20}%
\special{pa 2402 2564}%
\special{pa 2642 2564}%
\special{fp}%
\special{sh 1}%
\special{pa 2642 2564}%
\special{pa 2574 2544}%
\special{pa 2588 2564}%
\special{pa 2574 2584}%
\special{pa 2642 2564}%
\special{fp}%
%
\special{pn 20}%
\special{pa 2362 2682}%
\special{pa 2602 2682}%
\special{fp}%
\special{sh 1}%
\special{pa 2602 2682}%
\special{pa 2534 2662}%
\special{pa 2548 2682}%
\special{pa 2534 2702}%
\special{pa 2602 2682}%
\special{fp}%
%
\special{pn 20}%
\special{pa 2344 2768}%
\special{pa 2562 2768}%
\special{fp}%
\special{sh 1}%
\special{pa 2562 2768}%
\special{pa 2494 2748}%
\special{pa 2508 2768}%
\special{pa 2494 2788}%
\special{pa 2562 2768}%
\special{fp}%
%
\special{pn 8}%
\special{pa 2472 3738}%
\special{pa 2472 1448}%
\special{fp}%
\end{picture}%
\hspace{1truecm}}

\vspace{0.35truecm}

\centerline{{\bf Figure 4 $\,:\,$ Geodesic variation $\delta$}}

\vspace{0.5truecm}

\noindent
%
From $(3.9)$ and $(3.10)$, we obtain 
\begin{align*}
\left.\frac{\partial^2\rho^h}{\partial z_i\partial\overline z_j}\right\vert_{p}
=&(\nabla d\rho^h)_{p}^{\mathbb C}
\left(\frac{1}{2}(\ve_i^p-\sqrt{-1}{\bf J}\ve^p_i),\,\frac{1}{2}(\ve_j^{p}+\sqrt{-1}{\bf J}\ve_j^{p})\right)\\
=&\frac{1}{4}\left(\,(\nabla d\rho^h)_{p}((\ve^p_i,\ve_j^{p})+(\nabla d\rho^h)_{p}({\bf J}\ve^p_i,{\bf J}\ve_j^{p})\right.\\
&\left.-\sqrt{-1}(\nabla d\rho^h)_{p}({\bf J}\ve_i^{p},\ve_j^{p})+\sqrt{-1}(\nabla d\rho^h)_{p}(\ve_i^{p},{\bf J}\ve_j^{p})\,\right)\\
=&\left\{\begin{array}{l}
\displaystyle{
\begin{array}{r}
\displaystyle{\frac{1}{4}\,\left((\nabla^dd\rho^d)_{p}(\ve_i^{p},\ve_j^{p})+\tanh\lambda_a(\vZ)\cdot\lambda_a({\rm grad}\,\rho)_{\vZ})\,\delta_{ij}\right)\qquad}\\
((i,j)\in\{m_{a-1}+1,\cdots,m_a\}\,\,\,(a=1,\cdots,l))
\end{array}}\\
\displaystyle{\frac{1}{4}\,\left((\nabla^dd\rho^d)_{p}(\ve_i^{p},\ve_j^{p}\right)} \hspace{2truecm} ((i,j):\,\,{\rm other}).
\end{array}\right.
\end{align*}
\qed

%
%
%

\vspace{0.5truecm}

Also, we can show the following fact.  

\vspace{0.5truecm}

\noindent
{\bf Lemma 3.2.} {\sl Assume that $p\in(G^d/K)_{\rm reg}$ (i.e., $\vZ\in\mathcal C$).  Then, for any $\rho\in{\rm Conv}_W^+(D)$, we have 
$$(\nabla^dd\rho^d)_p((\ve^0_i)^p,(\ve^0_j)^p)=(\nabla^0d\rho)_{\vZ}(\ve^0_i,\ve^0_j)
=\left.\frac{\partial^2\rho}{\partial x_i\partial x_j}\right\vert_Z,\leqno{(3.11)}$$
where $(x_1,\cdots,x_r)$ is the Euclidean coordinate of $\mathfrak a^d$ with respect to $(\ve_1,\cdots,$\newline
$\ve_n)$ and 
$$(\nabla^dd\rho^d)_p((\ve^{\lambda_a}_i)^p,(\ve^{\lambda_b}_j)^p)=\frac{-1}{\tanh\lambda_a(\vZ)}\cdot\lambda_a(({\rm grad}\,\rho)_{\vZ})
\delta_{ab}\delta_{ij}.
\leqno{(3.12)}$$
Also, we have 
$$(\nabla^dd\rho^d)_p((\ve^{\lambda_a}_i)^p,(\ve^0_j)^p)=0.\leqno{(3.13)}$$
}

\vspace{0.5truecm}

\noindent
{\it Proof.}\ \ 
Since $p\in(G^d/K)_{\rm reg}$, the normal space $T_p^{\perp}{\rm Exp}_o(\mathfrak a^d)$ of the submanifold ${\rm Exp}_o(\mathfrak a^d)$ in $G^d/K$ 
at $p$ is equal to $T_p(K\cdot p)$.  Let $\gamma_i$ ($i=1,\cdots,r$) be the geodesic in ${\rm Exp}_o(\mathfrak a^d)$ with $\gamma_i'(0)=(\ve_i^0)^p$.  
Then, since ${\gamma}_i$ ($i=1,\cdots,r$) is a geodesic in $G^d/K$, we have 
\begin{align*}
(\nabla^dd\rho^d)_p((\ve_i^0)^p,(\ve_i^0)^p)&=\left.\frac{d^2(\rho^d\circ{\gamma}_i)}{ds^2}\right\vert_{s=0}-d\rho^d_p(\nabla^d_{\gamma_i'(0)}\gamma_i')\\
&=\left.\frac{d^2(\rho^d\circ{\gamma}_i)}{ds^2}\right\vert_{s=0}-d\rho_{\vZ}(\nabla^0_{\gamma_i'(0)}\gamma_i')\\
&=(\nabla^0d\rho)_{\vZ}(\ve_i^0,\ve_i^0).
\end{align*}
Similarly, we can show that 
$$(\nabla^dd\rho^d)_p((\ve_i^0)^p+(\ve_j^0)^p,(\ve_i^0)^p+(\ve_j^0)^p)=(\nabla^0d\rho)_{\vZ}(\ve_i^0+\ve_j^0,\ve_i^0+\ve_j^0).$$
Hence, from the symmetricnesses of $(\nabla^dd\rho^d)_p$ and $(\nabla^0d\rho)_{\vZ}$, we have 
$$(\nabla^dd\rho^d)_p((\ve_i^0)^p,(\ve_j^0)^p)=(\nabla^0d\rho)_Z(\ve_i^0,\ve_j^0).$$
Thus we obtain the desired relation $(3.11)$.  

Take $\vw\in T_p(K\cdot p)$ and let $\widehat{\gamma}$ be the geodesic in $K\cdot p$ with $\widehat{\gamma}'(0)=\vw$.  Then, 
since $\rho^d$ is $K$-invariant, we have 
\begin{align*}
(\nabla^dd\rho^d)_p(\vw,\vw)=\left.\frac{d^2(\rho^d\circ\widehat{\gamma})}{ds^2}\right\vert_{s=0}
-d\rho^d_p(h_p(\vw,\vw))=-d\rho_{\vZ}(h_p(\vw,\vw)).
\end{align*}
Hence, from the symmetricnesses of $(\nabla^dd\rho^d)_p$ and $d\rho_{\vZ}(h_p(\cdot,\cdot))$, we have 
$$(\nabla^dd\rho^d)_p(\vw_1,\vw_2)=-d\rho_{\vZ}(h_p(\vw_1,\vw_2))$$
for any $\vw_1,\vw_2\in T_p(K\cdot p)$.  
On the other hand, according to $(3.2)$, we have 
$$h_p((\ve^{\lambda_a}_i)^p,(\ve^{\lambda_b}_j)^p)=\delta_{ij}\delta_{ab}\cdot\frac{1}{\tanh\lambda_a(\vZ)}\,\lambda_a^{\sharp},\leqno{(3.14)}$$
where we used also $\beta_A((\ve^{\lambda_a}_i)^p,(\ve^{\lambda_b}_j)^p)=-\delta_{ab}\delta_{ij}$.   
Hence we obtain the desired relation $(3.12)$.  

Let $\gamma_i$ be the geodesic in $K\cdot p$ with $\gamma_i'(0)=(\ve^{\lambda_a}_i)^p$.  
This geodesic $\widehat{\gamma}_i$ is given as $\widehat{\gamma}_i(t):=\exp(tw)(p)$ for some $w\in\mathfrak k$.  
Define a vector field $\widetilde{\ve_j^0}$ along $\widehat{\gamma}_i$ by $\widetilde{\ve_j^0}(t):=\exp(tw)_{\ast}(\ve_j^0)$.  
It is shown that $\widetilde{\ve_j^0}$ is a $K$-equivariant $\nabla^{\perp}$-parallel normal vector field of $K\cdot p$ along $\gamma_i$ because $K$-action on $G^d/K$ 
is hyperpolar (see \cite{BS}, \cite{HPTT}), where $\nabla^{\perp}$ denotes the normal connection of the submanifold $K\cdot p$ in $G^d/K$.  
Since $\rho^d$ is $K$-invariant and $\widetilde{\ve_j^0}$ is $\nabla^{\perp}$-parallel, we have 
\begin{align*}
(\nabla^dd\rho^d)_p((\ve_i^{\lambda_a})^p,(\ve_j^0)^p)&=\left.\frac{d\widetilde{\ve_j^0}(t)(\rho^d)}{dt}\right\vert_{t=0}-(\nabla^d_{\gamma'_i(0)}\widetilde{\ve_j^0})(\rho^d)\\
&=((A_p)_{(\ve_j^0)^p}((\ve_i^{\lambda_a})^p)(\rho^d)=0.
\end{align*}
\qed
%

\vspace{0.5truecm}

Define a non-linear differential operator 
$\mathcal D:{\rm Conv}_W^+(\mathfrak a^d)\to C^{\infty}(W\cdot\mathcal C)$ of order one by 
$$\mathcal D(\rho):=(-1)^{n-r}\,\mathop{\Pi}_{\lambda\in\triangle_+}\left(\frac{2(\lambda\circ{\rm grad}\,\rho)\vert_{W\cdot\mathcal C}}
{\sinh\circ 2(\lambda\vert_{W\cdot\mathcal C})}\right)^{m_{\lambda}}
\quad(\rho\in{\rm Conv}_W^+(\mathfrak a^d)),\leqno{(3.15)}$$
where ${\rm grad}\,\rho$ is the gradient vector field of $\rho$ with respect to $\beta_0$.  
If $\rho\in{\rm Conv}_W^+(\mathfrak a^d)$, then ${\rm grad}\,\rho$ is $W$-equivariant.  This together with the $W$-invariance of $\lambda$'s implies that 
$\mathcal D(\rho)$ is $W$-invariant.  Hence we obtain $\mathcal D({\rm Conv}_W^+(\mathfrak a^d))\subset C^{\infty}_W(W\cdot\mathcal C)$.  

Let $\iota$ be the inclusion mapping of $G^d\cdot o(=G^d/K)$ into $G^{\mathbb C}/K^{\mathbb C}$.  
The structure of the K$\ddot{\rm a}$hler metric $\beta_{\rho^h}$ on $G^{\mathbb C}/K^{\mathbb C}$ is dominated by that of the induced metric $\iota^{\ast}\beta_{\rho^h}$ 
on $G\cdot o(=G^d/K)$.  
In fact, the anti-Kaehler metric $\beta_{\rho^h}$ is the real part of the holomorphic extension $(\iota^{\ast}\beta_{\rho^h})^h$ of 
the induced metric $\iota^{\ast}\beta_{\rho^h}$, where $(\iota^{\ast}\beta_{\rho^h})^h$ is the holomorphic metric (i.e., the holomorphic section of 
$((T(G^{\mathbb C}/K^{\mathbb C}))^{\mathbb C})^{\ast}\otimes((T(G^{\mathbb C}/K^{\mathbb C}))^{\mathbb C})^{\ast}$) such that the restriction of its real part 
${\rm Re}(\iota^{\ast}\beta_{\rho^h})^h\vert_{T^{\ast}(G^{\mathbb C}/K^{\mathbb C})\otimes T^{\ast}(G^{\mathbb C}/K^{\mathbb C})}$ to $G^d/K$ 
is equal to $\iota^{\ast}\beta_{\rho^h}$.  

From Lemmas 3.1 and 3.2, we obtain the following fact for this induced metric $\iota^{\ast}\beta_{\rho^h}$.  

\vspace{0.25truecm}

\noindent
{\bf Proposition 3.3.} {\sl {\rm (i)}\ \ For the induced metric $\iota^{\ast}\beta_{\rho^h}$, the following relation holds:
{\small
$$(\iota^{\ast}\beta_{\rho^h})_p(\ve_i^p,\ve_j^p)
=\left\{\begin{array}{ll}
\displaystyle{\frac{1}{2}\cdot\left.\frac{\partial^2\rho}{\partial x_i\partial x_j}\right\vert_Z} & (i,j\in\{1,\cdots,m_0\})\\
\displaystyle{-\frac{\lambda_1(({\rm grad}\,\rho)_{\vZ})}{\sinh 2\lambda_1(\vZ)}\cdot\delta_{ij}} & (i,j\in\{m_0+1,\cdots,m_1\})\\
\qquad\qquad\,\,\vdots & \qquad\qquad\qquad\vdots\\
\displaystyle{-\frac{\lambda_l(({\rm grad}\,\rho)_{\vZ})}{\sinh 2\lambda_l(\vZ)}\cdot\delta_{ij}} & (i,j\in\{m_{l-1}+1,\cdots,m_l\}).
\end{array}\right.
\leqno{(3.16)}$$
}

{\rm (ii)}\ \ For $p:={\rm Exp}_o(\vZ)$ ($\vZ\in\mathcal C$), the following relations hold:
$$
{\rm det}\left(\left.\frac{\partial^2\rho^h}{\partial z^p_i\partial\overline z^p_j}\right\vert_p\right)
=\frac{1}{4^n}\cdot{\rm det}\left(\left.\frac{\partial^2\rho}{\partial x_i\partial x_j}\right\vert_Z\right)\cdot\mathcal D(\rho)(\vZ),
\leqno{(3.17)}$$
where $(x_1,\cdots,x_r)$ is the Euclidean coordinate of $\mathfrak a^d$.
}

\vspace{0.25truecm}

\noindent
{\it Proof.}\ \ 
From the relations in Lemmas 3.1 and 3.2, we can derive 
{\small
$$\left.\frac{\partial^2\rho^h}{\partial z_i\partial\overline z_j}\right\vert_{p}
=\left\{\begin{array}{ll}
\displaystyle{\frac{1}{4}\cdot\left.\frac{\partial^2\rho}{\partial x_i\partial x_j}\right\vert_Z} & (i,j\in\{1,\cdots,m_0\})\\
\displaystyle{-\frac{\lambda_1(({\rm grad}\,\rho)_{\vZ})}{2\sinh 2\lambda_1(\vZ)}\cdot\delta_{ij}} & 
(i,j\in\{m_0+1,\cdots,m_1\})\\
\qquad\qquad\,\,\vdots & \qquad\qquad\qquad\vdots\\
\displaystyle{-\frac{\lambda_l(({\rm grad}\,\rho)_{\vZ})}{2\sinh 2\lambda_l(\vZ)}\cdot\delta_{ij}} & 
(i,j\in\{m_{l-1}+1,\cdots,m_l\}).
\end{array}\right.
\leqno{(3.18)}$$
}
From this relation, we can derive 
the relation in the statement (ii) directly.  
Also, by noticing that 
$$(\iota^{\ast}\beta_{\rho^h})_p(\ve_i^p,\ve_j^p)=\left.\frac{\partial^2\rho^h}{\partial z_i\partial\overline z_j}\right\vert_{p}
+\left.\frac{\partial^2\rho^h}{\partial z_j\partial\overline z_i}\right\vert_{p}$$
holds, we can derive the relation in the statement (i).  
\qed

\vspace{0.35truecm}

Let $(x_1,\cdots,x_r)$ be the Euclidean coordinate of $\mathfrak a^d$.  
From $(2.1)$ and $(3.17)$, we see that $\beta_{\rho^h}$ is Ricci-flat if the following relation 
$$\mathop{\Pi}_{\lambda\in\triangle_+}2^{m_{\lambda}}(\lambda\circ{\rm grad}\,\rho)^{m_{\lambda}}
\cdot{\rm det}\left(\frac{\partial^2\rho}{\partial x_i\partial x_j}\right)
=c\cdot\mathop{\Pi}_{\lambda\in\triangle_+}(\sinh\circ 2\lambda)^{m_{\lambda}}\leqno{(3.19)}$$
holds for some positive constant $c$.  
%
%

\vspace{0.25truecm}

\noindent
{\it Example 3.1.}\ \ We consider the case where $G^{\mathbb C}/K^{\mathbb C}=SO(n+1,\mathbb C)/SO(n,\mathbb C)$ ($n$-dimensional complex sphere).  
In this case, the restricted root system $\triangle$ of $G^d/K=SO_0(1,n)/SO(n)$ for the maximal abelian subspace $\mathfrak a^d:={\rm Span}\{\vZ\}$ of $\mathfrak p^d$ 
is given by 
$\triangle=\{\lambda\}$, where $\lambda$ is the element of $(\mathfrak a^d)^{\ast}$ defined by $\lambda(x\ve_1):=x$ ($x\in\mathbb R$) 
($\ve_1:=\frac{\vZ}{\vert\vert\vZ\vert\vert}$).  
Set $p:={\rm Exp}_o(\vZ)$ and let $\ve^p_i$ ($i=1,\cdots,n$) be as above.  Note that $\ve_i\in\mathfrak p_{\lambda}$ ($i=2,\cdots,n$).  
Hence $(3.19)$ is equivalent to 
$$2^{n-1}\widehat{\rho}'(x)^{n-1}\widehat{\rho}''(x)=c\,\sinh^{n-1}2x,\leqno{(3.20)}$$
where $\widehat{\rho}$ is a function over $\mathbb R$ defined by $\widehat{\rho}(x):=\rho(x\ve_1)$ ($x\in\mathbb R$).  

In particular, we shall consider the case of $n=2$.  In this case, the relation $(3.20)$ is as follows:
$$(\widehat{\rho}'(x)^2)'=c\sinh 2x.\leqno{(3.21)}$$
Hence we obtain $\widehat{\rho}(x)=\sqrt c\,\cosh x$.  
Note that, in this case, $\beta_{\rho^h}$ is the Eguchi-Hanson metric.  
According to $(3.16)$, the metric $\iota^{\ast}\beta_{\rho^h}$ on $G^d\cdot o$ induced from the Eguchi-Hanson metric is given by 
{\small
$$\left((\iota^{\ast}\beta_{\rho^h})_p(\ve_i^p,\ve_j^p)\right)
=\left(\begin{array}{cc}
\displaystyle{\frac{\sqrt{c}}{2}\,\cosh\vert\vert\vZ\vert\vert} & 0\\
0 & \displaystyle{-\frac{\sqrt{c}}{2\cosh\vert\vert\vZ\vert\vert}}
\end{array}\right),$$
}where $p$ and $\vZ$ are as above.  

\vspace{0.25truecm}

\noindent
{\it Example 3.2.}\ \ We consider the case where $G^{\mathbb C}/K^{\mathbb C}$ is one of rank one complex symmetric spaces 
$$\begin{array}{c}
SL(n+1,\mathbb C)/S(GL(1,\mathbb C)\times GL(n,\mathbb C)),\\
Sp(n+1,\mathbb C)/(Sp(1,\mathbb C)\times Sp(n,\mathbb C)),\,\,\,F_4^{\mathbb C}/{\rm SO}(9,\mathbb C)
\end{array}$$
other than complex spheres.  
In this case, the restricted root system $\triangle$ of 
\begin{align*}
G^d/K=&SU(1,n)/S(U(1)\times U(n-1)),\\
&Sp(1,n)/(Sp(1)\times Sp(n)\,\,{\rm or}\,\,F_4^{-20}/{\rm Spin}(9)
\end{align*}
for a maximal abelian subspace $\mathfrak a^d:={\rm Span}\{\vZ\}$ of $\mathfrak p^d$ is given by 
$\triangle=\{\lambda,2\lambda\}$ and 
$${\rm dim}\,\mathfrak p_{2\lambda}=\left\{\begin{array}{ll}
1 & ({\rm when}\,\,G^d/K=SU(1,n)/S(U(1)\times U(n)))\\
3 & ({\rm when}\,\,G^d/K=Sp(1,n)(Sp(1)\times Sp(n)))\\
7 & ({\rm when}\,\,G^d/K=F_4^{-20}/{\rm Spin}(9)),
\end{array}\right.$$
where $\lambda$ is the element of $(\mathfrak a^d)^{\ast}$ defined by $\lambda(x\ve_1):=x$ ($x\in\mathbb R$) ($\ve_1:=\frac{\vZ}{\vert\vert\vZ\vert\vert}$).  
For the convenience, we set 
$$d:=\left\{\begin{array}{ll}
1 & ({\rm when}\,\,G^d/K=SU(1,n)/S(U(1)\times U(n)))\\
3 & ({\rm when}\,\,G^d/K=Sp(1,n)(Sp(1)\times Sp(n)))\\
7 & ({\rm when}\,\,G^d/K=F_4^{-20}/{\rm Spin}(9)).
\end{array}\right.$$
Set $p:={\rm Exp}_o(\vZ)$ and let $\ve^p_i$ ($i=1,\cdots,(d+1)n$) be as above, where $n=2$ in the case of $G^{\mathbb C}/K^{\mathbb C}=F_4^{\mathbb C}/SO(9,\mathbb C)$.  
Note that $\ve_i\in\mathfrak p_{\lambda}$ ($i=2,\cdots,(d+1)n-d-1$) and $\ve_j\in\mathfrak p_{2\lambda}$ ($j=(d+1)n-d,\cdots,(d+1)n$).  
Hence $(3.19)$ is equivalent to 
$$2^{(d+1)n+d-1}\widehat{\rho}'(x)^{(d+1)n-1}\widehat{\rho}''(x)=c\,\sinh^{(d+1)n-d-1}(2x)\,\sinh^d(4x),$$
where $\widehat{\rho}$ is a function over $\mathbb R$ defined as above.  

\vspace{0.35truecm}

In 1990, L. A. Caffarelli proved the following fact (see Theorem 1-(b) and Theorem 2 in \cite{C3}). 

\vspace{0.35truecm}

\noindent
{\bf Proposition 3.4(\cite{C3}).} {\sl Let $D$ be a convex bounded domain of $\mathbb R^r$ with 
$B^r(\varepsilon)\subset D\subset B^r(r\varepsilon)$ ($\varepsilon>0$) and $f$ be a strictly positive $C^{0,\alpha}$-continuous function on $D$ ($\alpha>0$), 
where $B^r(\cdot)$ denotes the Euclidean ball of radius $(\cdot)$ centered at $(0,\cdots,0)$.  
If $\rho:D\to\mathbb R$ be a convex viscosity solution of 
$${\rm det}\left(\frac{\partial^2\rho}{\partial x_i\partial x_j}\right)=f\quad\,\,\,\,({\rm on}\,\,D)$$
($(x_1,\cdots,x_r)\,:\,$ the Euclidean coordinate of $\mathbb R^r$) satisfying $\rho\vert_{\partial D}=0$, then $\rho$ is of class $C^{2,\alpha}$ on 
$B^r(\frac{\varepsilon}{2})$.}

\vspace{0.35truecm}

In 2003, R. Bielawski proved the following fact (see Theorems A.1 and Corollary A.3 in \cite{B1}).  

\vspace{0.35truecm}

\noindent
{\bf Proposition 3.5(\cite{B1}).}
{\sl Let $W\curvearrowright\mathbb R^r$ be a finite irreducible reflection group and 
$F_1,F_2$ be non-negative locally bounded $W$-invariant measurable functions on $\mathbb R^r$.  
Assume that $\displaystyle{\int_{\mathbb R^r}F_1(x_1,\cdots,x_r)\,dx_1\cdots dx_r=\infty}$ holds.  
Then the following statements (i) and (ii) hold.  

(i)\ There exists a $W$-invariant global weak solution of the following Monge-Amp$\acute{\rm e}$re equation:
$$(F_1\circ{\rm grad}\,\rho)\cdot{\rm det}\left(\frac{\partial^2\rho}{\partial x_i\partial x_j}\right)=F_2,\leqno{(3.22)}$$
where the gradient vector field ${\rm grad}\,\rho$ is regarded as a (multi-valued) map of $\mathbb R^r$ to oneself.  
This weak solution $\rho$ is convex and Lipschitz continuous.  

(ii)\ If $F_1,F_2$ are positive functions of class $C^{k,\alpha}$ ($k\geq 0,\,\,\alpha>0$), then 
there exists a $W$-invariant strictly convex global $C^{k+2,\alpha}$-solution of $(3.22)$.  
Furthermore, this $W$-invariant global $C^{k+2,\alpha}$-solution is of class $C^{\infty}$ 
by reducing this problem to the fully non-linear uniformly elliptic case (\cite{C1}) and using the higher regularity (which is shown by the bootstrap argument) of 
the $C^{2,\alpha}$-solution of fully non-linear uniformly elliptic partial differential equation of order two.
}

\vspace{0.35truecm}

\noindent
{\it Remark 3.2.}\ 
(i)\ A {\it $W$-invariant convex global weak solution} of $(3.22)$ means a convex function $\rho$ on $\mathbb R^r$ such that 
$$\int_BF_2(x_1,\cdots,x_r)\,dx_1\cdots dx_r=\int_{({\rm grad}\,\rho)(B)}F_1(x_1,\cdots,x_r)\,dx_1\cdots dx_r\leqno{(3.23)}$$
holds for any Borel set $B$ of $\mathbb R^r$.  

(ii)\ The gradient vector field ${\rm grad}\,\rho$ is defined by 
$$\begin{array}{r}
\displaystyle{({\rm grad}\,\rho)_{\vx}:=\left\{\left.\left(\left(\frac{\partial\nu_P}{\partial x_1}\right)_{\vx},\cdots,\left(\frac{\partial\nu_P}{\partial x_r}\right)_{\vx}\right)\,
\right\vert\,P\in\mathcal S_{\vx}(\rho)\right\}}\\
(\vx=(x_1,\cdots,x_r)\in\mathbb R^r),
\end{array}$$
where $\mathcal S_{\vx}(\rho)$ is the set of all support hyperplane of the graph of $\rho$ at $\vx$ and $\nu_P$ is the affine function on $\mathbb R^r$ whose graph is equal to $P$.  
Here we note that $P$ is called a {\it support hyperplane of the graph of $\rho$ at $\vx$} if $\nu_P(\vx)=\rho(\vx)$ and $\nu_P\leq\rho$ hold.  

\vspace{0.35truecm}

Here we shall give some examples of the gradient vector in the above sense.  

\vspace{0.35truecm}

\noindent
{\it Example 3.3.}\ \ 
Let $\rho$ be a $C^{\infty}$-function on $\mathbb R^2$ defined by 
$$\rho(x_1,x_2):=x_1^2+x_2^2+1\quad\,\,((x_1,x_2)\in\mathbb R^2).$$
Since the graph of $\rho$ is given by $\{(x_1,x_2,y)\in\mathbb R^3\,\vert\,y=x_1^2+x_2^2+1\}$, we have 
$$\mathcal S_{(1,0)}(\rho):=\{P\}\quad\,\,(P:=\{(x_1,x_2,y)\in\mathbb R^3\,\vert\,y=2x_1\}).$$
Since $P$ is the graph of the affine function $\nu(x_1,x_2)=2x_1$, we obtain \newline
$({\rm grad}\,\rho)_{(1,0)}=\{(2,0)\}$.  

\vspace{0.35truecm}

\noindent
{\it Example 3.4.}\ \ 
Let $\rho$ be a continuous function on $\mathbb R^2$ defined by 
$$\rho(x_1,x_2):=\left\{\begin{array}{ll}
\sqrt{x_1^2+x_2^2}+1 & (0\leq x_1^2+x_2^2\leq 1)\\
2\sqrt{x_1^2+x_2^2} & (x_1^2+x_2^2\geq 1)
\end{array}\right.$$
Then we have 
$$\mathcal S_{(1,0)}(\rho)=\{P_a\,\vert\,1\leq a\leq 2\}\quad\,\,(P_a:=\{(x_1,x_2,y)\,\vert\,y=ax_1-a+2\}).$$
Since $P_a$ is the graph of the affine function $\nu_a(x_1,x_2)=ax_1-a+2$, we obtain 
$$({\rm grad}\,\rho)_{(1,0)}=\{(a,0)\,\vert\,1\leq a\leq 2\}.$$

\vspace{0.35truecm}

\noindent
{\it Example 3.5.}\ \ 
Let $\rho$ be a continuous function on $\mathbb R^2$ defined by 
$$\rho(x_1,x_2):=\left\{\begin{array}{ll}
\displaystyle{\vert x_1\vert+\vert x_2\vert+1} & ((0\leq\vert x_1\vert+\vert x_2\vert\leq 1)\\
\displaystyle{2(\vert x_1\vert+\vert x_2\vert)} & (\vert x_1\vert+\vert x_2\vert\geq 1)
\end{array}\right.$$
Then we have 
$$\mathcal S_{(1,0)}(\rho)=\left\{P_{a,\theta}\,\left\vert\,1\leq a\leq 2,\,\,-\frac{\pi}{4}\leq\theta\leq\frac{\pi}{4}\right.\right\},$$
where $P_{a,\theta}$ is the hyperplanes in $\mathbb R^3$ defined by 
$$P_{a,\theta}:=\{(x_1,x_2,y)\,\vert\,a\cos\theta\cdot(x_1-1)+a\sin\theta\cdot x_2-\cos\theta\cdot(y-2)=0\}.$$
Since $P_{a,\theta}$ is the graph of the affine function 
$\nu_{a,\theta}(x_1,x_2):=ax_1+a\tan\theta\cdot x_2+2-a$, we obtain 
\begin{align*}
({\rm grad}\,\rho)_{(1,0)}&=\left\{(a,a\tan\theta)\,\left\vert\,1\leq a\leq 2,\,\,-\frac{\pi}{4}\leq\theta\leq\frac{\pi}{4}\right.\right\}\\
&=\{(a,ab)\,\vert\,1\leq a\leq 2,\,\,-1\leq b\leq 1\}.
\end{align*}

In 1991, L. A. Caffarelli (\cite{C4}) proved the following fact for the regularity of locally Lipschitz functions on the whole of a Euclidean space satisfying 
some Monge-Amp$\grave{\rm e}$re inequality (see Theorem 2 of \cite{C4}).  

\vspace{0.35truecm}

\noindent
{\bf Proposition 3.6(\cite{C4}).}\ \ 
{\sl Let $\rho$ be a locally Lipschitz function on $\mathbb R^r$ satisfying the Monge-Amp$\grave{\rm e}$re inequality:
$$C_1\leq{\rm det}\left(\frac{\partial^2\rho}{\partial x_i\partial x_j}\right)\leq C_2(<\infty)$$
($C_1,C_2\,:\,$positive constants).  If $\rho$ is strictly convex in the usual sense, then it is locally of class $C^{1,\alpha}$, 
where ``locally'' means that $\alpha$ can be chosen uniformly only on a neighborhood of each point of $\mathbb R^r$.
}

\vspace{0.35truecm}

\noindent
{\it Remark 3.3.}\ 
The Monge-Amp$\grave{\rm e}$re inequality in the statement of Proposition 3.6 means that, 
$$C_1\cdot{\rm Vol}(B)\leq{\rm Vol}(({\rm grad}\,\rho)(B))\leq C_2\cdot{\rm Vol}(B)$$
holds for any Borel set $B$ of $\mathbb R^r$, where ${\rm Vol}(\cdot)$ denotes the Euclidean volume of $(\cdot)$.  

\vspace{0.35truecm}

Define functions $\widehat F_1$ and $\widehat F_2$ on $\mathfrak a^d$ by 
$$
\widehat F_1(\vZ):=\mathop{\Pi}_{\lambda\in\triangle_+}2^{m_{\lambda}}\vert\lambda(\vZ)\vert^{m_{\lambda}},\quad
\widehat F_2(\vZ):=c\cdot\mathop{\Pi}_{\lambda\in\triangle_+}\vert\sinh 2\lambda(\vZ))\vert^{m_{\lambda}}
$$
($\vZ\in\mathfrak a^d$).  Then the relation $(3.19)$ is rewritten as 
$$(\widehat F_1\circ{\rm grad}\,\rho)\cdot{\rm det}\left(\frac{\partial^2\rho}{\partial x_i\partial x_j}\right)=\widehat F_2.\leqno{(3.24)}$$
It is clear that $\widehat F_i$ ($i=1,2$) are locally bounded and non-negative $W$-invariant continuous (hence measurable) functions on $\mathfrak a^d$ and that $\widehat F_1$ 
satisfies\newline
$\displaystyle{\int_{\mathfrak a^d}\widehat F_1(x_1,\cdots,x_r)\,dx_1\cdots dx_r=\infty}$.  Hence we can apply the statement (i) of Proposition 3.5 to $(3.24)$ and 
show the existence of a $W$-invariant convex and Lipschitz continuous global weak solution of $(3.24)$.  
Regrettably, $\widehat F_i$ ($i=1,2$) are not positive.  In fact, they are equal to zero along $\mathfrak a^d\setminus W\cdot\mathcal C$.  
Also, they are not necessarily of class $C^1$ along $\mathfrak a^d\setminus W\cdot\mathcal C$.  
Hence we cannot apply the statement (ii) of Proposition 3.5 to $(3.24)$.  
However, we can derive the statement of Theorem A (stated in Introduction) in the case where $G/K$ is of rank two.  
We prove Theorem A by using Propositions 3.3, 3.4 and 3.5.  

\vspace{0.35truecm}

\noindent
{\it Proof of Theorem A.}\ As above, since 
$\widehat F_i$ ($i=1,2$) are locally bounded and non-negative $W$-invariant continuous (hence measurable) functions on $\mathfrak a^d$ and that $\widehat F_1$ satisfies \newline
$\displaystyle{\int_{\mathfrak a^d}\widehat F_1(x_1,x_2)\,dx_1dx_2=\infty}$, 
it follows from (i) of Proposition 3.5 that there exists a $W$-invariant 
weak solution $\rho$ of $(3.24)$ defined on the whole of $\mathfrak a^d$ and it is convex and Lipschitz continuous.  
Set $$\theta:=\left\{\begin{array}{ll}
\displaystyle{\frac{\pi}{3}} & ({\rm when}\,\,\triangle:(\mathfrak a_2){\rm -type})\\
\displaystyle{\frac{\pi}{4}} & ({\rm when}\,\,\triangle:(\mathfrak b_2){\rm -type})\\
\displaystyle{\frac{\pi}{6}} & ({\rm when}\,\,\triangle:(\mathfrak g_2){\rm -type})\\
\displaystyle{\frac{\pi}{2}} & ({\rm when}\,\,\triangle:(\mathfrak a_1\times\mathfrak a_1){\rm -type}).
\end{array}\right.$$
Also, let $\Pi=\{\lambda_1,\lambda_2\}\,(\subset\triangle_+)$ be the fundamental root system.  
We can take an orthonormal basis $(\ve_1,\ve_2)$ of $\mathfrak a^d$ such that $\lambda_1(\sum\limits_{i=1}^2x_i\ve_i)=x_2$ and 
$\lambda_2(\sum\limits_{i=1}^2x_i\ve_i)=x_1\sin\theta-x_2\cos\theta$ ($x_1,x_2\in\mathbb R$) hold.  
Then a Weyl domain $\mathcal C$ is given by 
\begin{align*}
\mathcal C&=\{\vZ\,\vert\,\lambda_i(\vZ)>0\,\,\,\,(i=1,2)\}\\
&=\left\{\left.\sum_{i=1}^2x_i\ve_i\,\right\vert\,x_2>0\,\,\,\&\,\,\,\frac{x_1}{x_2}>\frac{1}{\tan\theta}\right\},
\end{align*}
where $\frac{1}{\tan\theta}$ implies $0$ in the case of $\theta=\frac{\pi}{2}$.  
Since $\rho$ is a convex weak solution of $(3.24)$ and $\widehat F_2$ is positive over $\mathcal C$, 
$\widehat F_1$ also is positive on $({\rm grad}\,\rho)(\mathcal C)$ in the weak sense, that is, 
$$\int_{({\rm grad}\,\rho)(B)}\widehat F_1(x_1,x_2)\,dx_1dx_2>0$$
holds for any Borel set $B$ with positive measure of $\mathcal C$.  
From this fact, it follows that ${\rm Vol}(({\rm grad}\,\rho)(B))>0$ holds for any Borel set $B$ with positive measure of $W\cdot\mathcal C$.  
Under the correspondence $a\ve_1+b\ve_2\,\,\leftrightarrow\,\,a+\sqrt{-1}b$ ($a,b\in\mathbb R$), we identify $\mathfrak a^d$ with the complex number field $\mathbb C$.  
Let $\vP$ be the position vector field on $\mathfrak a^d$ (i.e., $\vP_{\vZ}:=\vZ$ ($\vZ\in\mathfrak a^d$)).  
Also, since $\rho$ is a $W$-invariant function and $W$ is generated by the reflections with respect to $\lambda^{-1}(0)$'s 
($\lambda\in\triangle_+$), ${\rm grad}\,\rho$ also is $W$-invariant and hence it is invariant under the reflections with respect to 
$\lambda^{-1}(0)$'s ($\lambda\in\triangle_+$).  Here the $W$-invariance of ${\rm grad}\,\rho$ means that 
$({\rm grad}\,\rho)_{w\cdot\vZ}=w\cdot({\rm grad}\,\rho)_{\vZ}$ holds for any $\vZ\in\mathcal C$ and any $w\in W$.  
In particular, if $w$ is the reflection with respect to $\lambda^{-1}(0)$ and if $\vZ\in\lambda^{-1}(0)$, then we have 
$({\rm grad}\,\rho)_{\vZ}=w\cdot({\rm grad}\,\rho)_{\vZ}$ holds and hence $\frac{({\rm grad}\,\rho)_{\vZ}}{\vert\vert({\rm grad}\,\rho)_{\vZ}\vert\vert}$ is described as 
$$\frac{({\rm grad}\,\rho)_{\vZ}}{\vert\vert({\rm grad}\,\rho)_{\vZ}\vert\vert}=\left\{\left.e^{\sqrt{-1}\theta}\frac{P_{\vZ}}{\vert\vert\vZ\vert\vert}\,\,\right\vert\,\,
-\varepsilon_{\vZ}\leq\theta\leq\varepsilon_{\vZ}\right\},$$
where $\varepsilon_{\vZ}$ is a non-negative number.  It is easy to show that $\varepsilon_{\vZ}=0$ for any $\vZ\in\partial\mathcal C$.  
On the other hand, since $\rho$ is convex, for any distinct 
$\vZ_1=e^{\sqrt{-1}\theta_1},\vZ_2=e^{\sqrt{-1}\theta_2}\in\mathfrak a^d$ ($0\leq\theta_1<\theta_2<2\pi$), $\frac{({\rm grad}\,\rho)_{\vZ_1}}{\vert\vert({\rm grad}\,\rho)_{\vZ_1}\vert\vert}$ 
and $\frac{({\rm grad}\,\rho)_{\vZ_1}}{\vert\vert({\rm grad}\,\rho)_{\vZ_2}\vert\vert}$ do not intersect.  
In more detail, they are described as 
$$\begin{array}{r}
\displaystyle{\frac{({\rm grad}\,\rho)_{\vZ_i}}{\vert\vert({\rm grad}\,\rho)_{\vZ_i}\vert\vert}=\left\{\left.e^{\sqrt{-1}\theta}\frac{P_{\vZ_i}}{\vert\vert\vZ_i\vert\vert}\,\,\right\vert\,\,-\varepsilon_i^-\leq\theta\leq
\varepsilon_i^+\right\}\quad\,\,(i=1,2)}\\
\displaystyle{(\theta_1+\varepsilon_1^+<\theta_2-\varepsilon_2^-),}
\end{array}$$
where $\varepsilon_i^{\pm}$ ($i=1,2$) are non-negative constants.  
Again, since $\rho$ is convex, we see that $\vert\vert({\rm grad}\,\rho)_{\vZ}\vert\vert$ is equal to some closed interval $[(a_1)_{\vZ},(a_2)_{\vZ}]$ ($(a_1)_{\vZ}>0$) for any 
$\vZ\in\mathfrak a^d\setminus\{{\bf 0}\}$.  
From these facts, we can derive that $({\rm grad}\,\rho)(w\cdot\mathcal C)\subset w\cdot\mathcal C$ holds for any $w\in W$ 
(see Figure 5 (when $\triangle$ is of type ($\mathfrak a_2$)), 
where we note that 
$$\mathfrak a^d\setminus W\cdot\mathcal C=\mathop{\cup}_{\lambda\in\triangle_+}\lambda^{-1}(0).$$
Hence, since $\lambda\circ{\rm grad}\,\rho\not=0$ holds on $W\cdot\mathcal C$ for any $\lambda\in\triangle_+$, 
$\widehat F_1\circ{\rm grad}\,\rho\not=0$ holds on $W\cdot\mathcal C$.  Hence, on $W\cdot\mathcal C$, $(3.24)$ is rewritten as 
$$\begin{array}{l}
\displaystyle{{\rm det}\left(\frac{\partial^2\rho}{\partial x_i\partial x_j}\right)
=\frac{\widehat F_2}{\widehat F_1\circ{\rm grad}\,\rho}}\\
\hspace{2.315truecm}\displaystyle{
=c\cdot\mathop{\Pi}_{\lambda\in\triangle_+}\left(\frac{\sinh\circ 2\lambda}{2(\lambda\circ{\rm grad}\,\rho)}\right)^{m_{\lambda}}
=\frac{c}{\mathcal D(\rho)}.}
\end{array}\leqno{(3.25)}$$
%
%
%
%
Note that 
$$\partial\mathcal C=\mathop{\amalg}_{i=1}^2(\lambda_i^{-1}(0)\cap\partial\mathcal C)\,\,\,\,\,({\rm hence}\,\,\,\,
(\mathfrak a^d\setminus W\cdot\mathcal C)\setminus\{0\}=\mathop{\amalg}_{i=1}^2W\cdot(\lambda_i^{-1}(0)\cap\partial\mathcal C))$$
holds.  For the simplicity, set $(\partial\mathcal C)_i:=\lambda_i^{-1}(0)\cap\partial\mathcal C$ ($i=1,2$).  
Since $\rho$ is a $W$-invariant convex function on $\mathfrak a^d$ and 
$\displaystyle{\int_{\mathfrak a^d}\widehat F_2(x_1,x_2)}\,dx_1dx_2>0$, it is shown that $\rho$ is a proper function (see Proposition 1.2 of \cite{B2}).  
On the basis of $(3.25)$ and $\varepsilon_{\vZ}={\bf 0}$ ($\vZ\in\partial\mathcal C$),
the natural extension of $\displaystyle{\frac{\widehat F_2}{\widehat F_1\circ{\rm grad}\,\rho}}$ on $W\cdot\mathcal C$ to $\mathfrak a^d$ is defined 
as a multi-valued function $\widehat F^e$ over $\mathfrak a^d$ given by 
$$\widehat F^e(w\cdot\vZ)
=\left\{\begin{array}{ll}
\displaystyle{\frac{\widehat F_2(\vZ)}{\widehat F_1(({\rm grad}\,\rho)_{\vZ})}} & (\vZ\in\mathcal C)\\
\displaystyle{c\cdot\left(\frac{\vert\vert\vZ\vert\vert}{\vert\vert({\rm grad}\,\rho)_{\vZ}\vert\vert}\right)^{m_{\lambda_i}+m_{2\lambda_i}}} & \\
\displaystyle{\times\mathop{\Pi}_{\lambda\in\triangle_+\setminus\{\lambda_i,2\lambda_i\}}
\left(\frac{\sinh(2\lambda(\vZ))}{2\lambda(({\rm grad}\,\rho)_{\vZ})}\right)^{m_{\lambda}}} & (\vZ\in(\partial\mathcal C)_i),\\
\displaystyle{c\cdot\left(\lim_{\vZ\to{\bf 0}}\frac{\vert\vert\vZ\vert\vert}{\vert\vert({\rm grad}\,\rho)_{\vZ}\vert\vert}\right)^{n-2}} & \\
\displaystyle{\times\mathop{\Pi}_{\lambda\in\triangle_+}\left(\lim_{\vZ\to{\bf 0}}
\frac{\sinh\left(2\vert\vert\vZ\vert\vert\lambda\left(\frac{\vZ}{\vert\vert\vZ\vert\vert}\right)\right)}
{2\vert\vert\vZ\vert\vert\lambda\left(\frac{({\rm grad}\,\rho)_{\vZ}}{\vert\vert({\rm grad}\,\rho)_{\vZ}\vert\vert}\right)}\right)^{m_{\lambda}}} & (\vZ={\bf 0}),
\end{array}\right.$$
($w\in W$), where $m_{2\lambda_i}=0$ in the case of $2\lambda_i\notin\triangle_+$.  

\vspace{1truecm}

\centerline{
\unitlength 0.1in
\begin{picture}(122.2700, 28.7000)(-71.2000,-36.0000)
%
\special{pn 8}%
\special{pa 2550 2256}%
\special{pa 4404 2256}%
\special{fp}%
\special{sh 1}%
\special{pa 4404 2256}%
\special{pa 4338 2236}%
\special{pa 4352 2256}%
\special{pa 4338 2276}%
\special{pa 4404 2256}%
\special{fp}%
%
\special{pn 8}%
\special{pa 3390 3062}%
\special{pa 3390 1350}%
\special{fp}%
\special{sh 1}%
\special{pa 3390 1350}%
\special{pa 3370 1416}%
\special{pa 3390 1402}%
\special{pa 3410 1416}%
\special{pa 3390 1350}%
\special{fp}%
%
\special{pn 8}%
\special{pa 3934 3046}%
\special{pa 2768 1342}%
\special{fp}%
%
\special{pn 20}%
\special{pa 3666 1858}%
\special{pa 3792 1658}%
\special{fp}%
\special{sh 1}%
\special{pa 3792 1658}%
\special{pa 3740 1704}%
\special{pa 3764 1702}%
\special{pa 3774 1724}%
\special{pa 3792 1658}%
\special{fp}%
%
\special{pn 20}%
\special{pa 3892 2256}%
\special{pa 4170 2256}%
\special{fp}%
\special{sh 1}%
\special{pa 4170 2256}%
\special{pa 4102 2236}%
\special{pa 4116 2256}%
\special{pa 4102 2276}%
\special{pa 4170 2256}%
\special{fp}%
%
\special{pn 20}%
\special{pa 3888 2126}%
\special{pa 4142 2084}%
\special{fp}%
\special{sh 1}%
\special{pa 4142 2084}%
\special{pa 4072 2076}%
\special{pa 4088 2092}%
\special{pa 4078 2114}%
\special{pa 4142 2084}%
\special{fp}%
\put(34.7800,-13.3200){\makebox(0,0)[rb]{$x_2$}}%
\put(44.6800,-21.8900){\makebox(0,0)[lt]{$x_1$}}%
%
\special{pn 8}%
\special{ar 4020 1390 150 146  2.5738149 4.5812442}%
%
\special{pn 8}%
\special{pa 3902 1462}%
\special{pa 3918 1480}%
\special{fp}%
\special{sh 1}%
\special{pa 3918 1480}%
\special{pa 3886 1418}%
\special{pa 3880 1442}%
\special{pa 3858 1446}%
\special{pa 3918 1480}%
\special{fp}%
%
\special{pn 8}%
\special{ar 4504 2214 228 178  3.1748030 4.6989023}%
%
\special{pn 8}%
\special{pa 4270 2214}%
\special{pa 4262 2246}%
\special{fp}%
\special{sh 1}%
\special{pa 4262 2246}%
\special{pa 4298 2188}%
\special{pa 4276 2196}%
\special{pa 4260 2176}%
\special{pa 4262 2246}%
\special{fp}%
\put(45.2100,-21.0900){\makebox(0,0)[lb]{$\lambda_1^{-1}(0)$}}%
\put(40.4300,-12.8400){\makebox(0,0)[lb]{$\lambda_2^{-1}(0)$}}%
%
\special{pn 20}%
\special{pa 2850 3382}%
\special{pa 2968 3220}%
\special{fp}%
\special{sh 1}%
\special{pa 2968 3220}%
\special{pa 2912 3262}%
\special{pa 2936 3264}%
\special{pa 2946 3286}%
\special{pa 2968 3220}%
\special{fp}%
\put(30.1300,-32.2000){\makebox(0,0)[lt]{'s $\displaystyle{,\,\,\left.\frac{{\rm grad}\,\rho}{\vert\vert{\rm grad}\,\rho\vert\vert}\right\vert_{\rho^{-1}(b)}}$}}%
%
\special{pn 4}%
\special{pa 3440 2256}%
\special{pa 3410 2226}%
\special{dt 0.027}%
\special{pa 3490 2256}%
\special{pa 3430 2198}%
\special{dt 0.027}%
\special{pa 3540 2256}%
\special{pa 3450 2166}%
\special{dt 0.027}%
\special{pa 3590 2256}%
\special{pa 3470 2138}%
\special{dt 0.027}%
\special{pa 3642 2256}%
\special{pa 3490 2110}%
\special{dt 0.027}%
\special{pa 3690 2256}%
\special{pa 3510 2080}%
\special{dt 0.027}%
\special{pa 3742 2256}%
\special{pa 3530 2052}%
\special{dt 0.027}%
\special{pa 3792 2256}%
\special{pa 3550 2022}%
\special{dt 0.027}%
\special{pa 3842 2256}%
\special{pa 3570 1992}%
\special{dt 0.027}%
\special{pa 3892 2256}%
\special{pa 3590 1964}%
\special{dt 0.027}%
\special{pa 3942 2256}%
\special{pa 3610 1934}%
\special{dt 0.027}%
\special{pa 3994 2256}%
\special{pa 3632 1906}%
\special{dt 0.027}%
\special{pa 4044 2256}%
\special{pa 3652 1876}%
\special{dt 0.027}%
\special{pa 4094 2256}%
\special{pa 3672 1846}%
\special{dt 0.027}%
\special{pa 4144 2256}%
\special{pa 3690 1818}%
\special{dt 0.027}%
\special{pa 4194 2256}%
\special{pa 3712 1790}%
\special{dt 0.027}%
\special{pa 4244 2256}%
\special{pa 3732 1760}%
\special{dt 0.027}%
\special{pa 4296 2256}%
\special{pa 3752 1730}%
\special{dt 0.027}%
\special{pa 4346 2256}%
\special{pa 3772 1702}%
\special{dt 0.027}%
\special{pa 4396 2256}%
\special{pa 3792 1672}%
\special{dt 0.027}%
\special{pa 4404 2214}%
\special{pa 3812 1644}%
\special{dt 0.027}%
\special{pa 4404 2166}%
\special{pa 3832 1614}%
\special{dt 0.027}%
\special{pa 4404 2116}%
\special{pa 3852 1586}%
\special{dt 0.027}%
\special{pa 4404 2068}%
\special{pa 3874 1556}%
\special{dt 0.027}%
\special{pa 4404 2022}%
\special{pa 3892 1528}%
\special{dt 0.027}%
\special{pa 4404 1972}%
\special{pa 3912 1498}%
\special{dt 0.027}%
\special{pa 4404 1924}%
\special{pa 3932 1470}%
\special{dt 0.027}%
\special{pa 4404 1876}%
\special{pa 3954 1440}%
\special{dt 0.027}%
\special{pa 4404 1826}%
\special{pa 3974 1410}%
\special{dt 0.027}%
\special{pa 4404 1778}%
\special{pa 3994 1382}%
\special{dt 0.027}%
%
\special{pn 4}%
\special{pa 4404 1730}%
\special{pa 4020 1358}%
\special{dt 0.027}%
\special{pa 4404 1682}%
\special{pa 4068 1358}%
\special{dt 0.027}%
\special{pa 4404 1632}%
\special{pa 4120 1358}%
\special{dt 0.027}%
\special{pa 4404 1584}%
\special{pa 4170 1358}%
\special{dt 0.027}%
\special{pa 4404 1536}%
\special{pa 4220 1358}%
\special{dt 0.027}%
\special{pa 4404 1486}%
\special{pa 4270 1358}%
\special{dt 0.027}%
\special{pa 4404 1438}%
\special{pa 4320 1358}%
\special{dt 0.027}%
\special{pa 4404 1390}%
\special{pa 4370 1358}%
\special{dt 0.027}%
\put(45.8200,-13.5000){\makebox(0,0)[lt]{$\mathcal C$}}%
%
\special{pn 8}%
\special{ar 2752 1270 220 284  1.0757581 3.0339028}%
%
\special{pn 8}%
\special{pa 2852 1512}%
\special{pa 2868 1486}%
\special{fp}%
\special{sh 1}%
\special{pa 2868 1486}%
\special{pa 2814 1530}%
\special{pa 2838 1530}%
\special{pa 2848 1552}%
\special{pa 2868 1486}%
\special{fp}%
\put(27.6500,-11.4900){\makebox(0,0)[lb]{$(\lambda_1+\lambda_2)^{-1}(0)$}}%
\put(43.6300,-28.9200){\makebox(0,0)[lb]{$\rho^{-1}(b)$}}%
%
\special{pn 8}%
\special{ar 3406 2256 488 622  5.8829058 6.2831853}%
%
\special{pn 8}%
\special{pa 3658 1858}%
\special{pa 3684 1876}%
\special{pa 3712 1892}%
\special{pa 3738 1910}%
\special{pa 3762 1928}%
\special{pa 3786 1950}%
\special{pa 3810 1970}%
\special{pa 3834 1990}%
\special{pa 3854 2014}%
\special{sp}%
%
\special{pn 8}%
\special{pa 3658 1858}%
\special{pa 3628 1844}%
\special{pa 3600 1832}%
\special{pa 3570 1820}%
\special{pa 3540 1808}%
\special{pa 3510 1798}%
\special{pa 3480 1788}%
\special{pa 3448 1784}%
\special{pa 3418 1774}%
\special{pa 3388 1772}%
\special{pa 3386 1772}%
\special{sp}%
%
\special{pn 8}%
\special{pa 3120 1866}%
\special{pa 3148 1852}%
\special{pa 3176 1838}%
\special{pa 3206 1826}%
\special{pa 3234 1812}%
\special{pa 3264 1802}%
\special{pa 3294 1792}%
\special{pa 3326 1784}%
\special{pa 3356 1776}%
\special{pa 3388 1772}%
\special{pa 3388 1772}%
\special{sp}%
%
\special{pn 8}%
\special{ar 3406 2270 540 644  3.1415927 3.5404996}%
%
\special{pn 8}%
\special{pa 3112 1866}%
\special{pa 3086 1882}%
\special{pa 3058 1900}%
\special{pa 3032 1918}%
\special{pa 3008 1938}%
\special{pa 2984 1958}%
\special{pa 2960 1978}%
\special{pa 2936 2000}%
\special{pa 2916 2022}%
\special{sp}%
%
\special{pn 8}%
\special{ar 3406 2262 488 622  6.2831853 6.2831853}%
\special{ar 3406 2262 488 622  0.0000000 0.4002795}%
%
\special{pn 8}%
\special{pa 3658 2650}%
\special{pa 3684 2634}%
\special{pa 3712 2616}%
\special{pa 3738 2598}%
\special{pa 3762 2580}%
\special{pa 3788 2562}%
\special{pa 3810 2540}%
\special{pa 3834 2518}%
\special{pa 3854 2496}%
\special{sp}%
%
\special{pn 8}%
\special{pa 3146 2644}%
\special{pa 3174 2658}%
\special{pa 3204 2672}%
\special{pa 3232 2686}%
\special{pa 3262 2698}%
\special{pa 3292 2708}%
\special{pa 3322 2718}%
\special{pa 3352 2726}%
\special{pa 3384 2736}%
\special{pa 3406 2740}%
\special{sp}%
%
\special{pn 8}%
\special{pa 3650 2654}%
\special{pa 3620 2668}%
\special{pa 3592 2680}%
\special{pa 3562 2692}%
\special{pa 3532 2704}%
\special{pa 3502 2714}%
\special{pa 3472 2722}%
\special{pa 3440 2730}%
\special{pa 3410 2738}%
\special{pa 3388 2738}%
\special{sp}%
%
\special{pn 8}%
\special{ar 3738 2246 872 666  2.7421776 3.1415927}%
%
\special{pn 8}%
\special{pa 3128 2634}%
\special{pa 3100 2620}%
\special{pa 3072 2604}%
\special{pa 3046 2588}%
\special{pa 3020 2570}%
\special{pa 2994 2550}%
\special{pa 2970 2530}%
\special{pa 2946 2510}%
\special{pa 2942 2504}%
\special{sp}%
%
\special{pn 20}%
\special{pa 3658 2650}%
\special{pa 3796 2846}%
\special{fp}%
\special{sh 1}%
\special{pa 3796 2846}%
\special{pa 3774 2780}%
\special{pa 3764 2802}%
\special{pa 3740 2802}%
\special{pa 3796 2846}%
\special{fp}%
%
\special{pn 20}%
\special{pa 3816 2562}%
\special{pa 3968 2748}%
\special{fp}%
\special{sh 1}%
\special{pa 3968 2748}%
\special{pa 3940 2684}%
\special{pa 3934 2706}%
\special{pa 3910 2708}%
\special{pa 3968 2748}%
\special{fp}%
%
\special{pn 20}%
\special{pa 3888 2372}%
\special{pa 4144 2390}%
\special{fp}%
\special{sh 1}%
\special{pa 4144 2390}%
\special{pa 4080 2366}%
\special{pa 4092 2386}%
\special{pa 4076 2406}%
\special{pa 4144 2390}%
\special{fp}%
%
\special{pn 20}%
\special{pa 3122 1866}%
\special{pa 2996 1664}%
\special{fp}%
\special{sh 1}%
\special{pa 2996 1664}%
\special{pa 3014 1732}%
\special{pa 3024 1710}%
\special{pa 3048 1710}%
\special{pa 2996 1664}%
\special{fp}%
%
\special{pn 20}%
\special{pa 2878 2256}%
\special{pa 2600 2256}%
\special{fp}%
\special{sh 1}%
\special{pa 2600 2256}%
\special{pa 2668 2276}%
\special{pa 2654 2256}%
\special{pa 2668 2236}%
\special{pa 2600 2256}%
\special{fp}%
\put(42.2000,-9.0000){\makebox(0,0)[rb]{{\small The existence range of $\displaystyle{\frac{({\rm grad}\,\rho)_{Z_0}}{\vert\vert({\rm grad}\,\rho)_{Z_0}\vert\vert}}$}}}%
\put(38.1100,-20.0700){\makebox(0,0)[rt]{{\small $\vZ_0$}}}%
%
\special{pn 8}%
\special{ar 2756 1260 230 192  2.8800949 4.5867108}%
%
\special{pn 8}%
\special{ar 4554 1888 400 492  3.8033928 4.6957516}%
%
\special{pn 8}%
\special{ar 4150 1380 958 564  4.8330383 6.2831853}%
\special{ar 4150 1380 958 564  0.0000000 1.5707963}%
%
\special{pn 8}%
\special{pa 4240 1588}%
\special{pa 4212 1624}%
\special{fp}%
\special{sh 1}%
\special{pa 4212 1624}%
\special{pa 4268 1584}%
\special{pa 4244 1582}%
\special{pa 4236 1560}%
\special{pa 4212 1624}%
\special{fp}%
%
\special{pn 20}%
\special{pa 3764 1916}%
\special{pa 3914 1730}%
\special{fp}%
\special{sh 1}%
\special{pa 3914 1730}%
\special{pa 3856 1770}%
\special{pa 3880 1772}%
\special{pa 3888 1794}%
\special{pa 3914 1730}%
\special{fp}%
%
\special{pn 20}%
\special{pa 3878 2010}%
\special{pa 4082 1878}%
\special{fp}%
\special{sh 1}%
\special{pa 4082 1878}%
\special{pa 4014 1898}%
\special{pa 4036 1908}%
\special{pa 4036 1932}%
\special{pa 4082 1878}%
\special{fp}%
%
\special{pn 20}%
\special{ar 3860 2020 276 256  5.3585079 6.0642517}%
%
\special{pn 8}%
\special{pa 3878 2010}%
\special{pa 4126 1956}%
\special{fp}%
%
\special{pn 8}%
\special{pa 3868 2000}%
\special{pa 4040 1818}%
\special{fp}%
%
\special{pn 20}%
\special{sh 1}%
\special{ar 3860 2016 10 10 0  6.28318530717959E+0000}%
\special{sh 1}%
\special{ar 3860 2016 10 10 0  6.28318530717959E+0000}%
\special{sh 1}%
\special{ar 3860 2016 10 10 0  6.28318530717959E+0000}%
%
\special{pn 20}%
\special{pa 3388 2730}%
\special{pa 3388 2966}%
\special{fp}%
\special{sh 1}%
\special{pa 3388 2966}%
\special{pa 3408 2898}%
\special{pa 3388 2912}%
\special{pa 3368 2898}%
\special{pa 3388 2966}%
\special{fp}%
%
\special{pn 20}%
\special{pa 3326 2956}%
\special{pa 3354 2966}%
\special{pa 3386 2974}%
\special{pa 3418 2974}%
\special{pa 3448 2968}%
\special{pa 3478 2958}%
\special{pa 3478 2956}%
\special{sp}%
%
\special{pn 8}%
\special{pa 3384 2728}%
\special{pa 3316 2962}%
\special{fp}%
%
\special{pn 8}%
\special{pa 3402 2728}%
\special{pa 3482 2960}%
\special{fp}%
%
\special{pn 20}%
\special{pa 2924 2000}%
\special{pa 2704 1894}%
\special{fp}%
\special{sh 1}%
\special{pa 2704 1894}%
\special{pa 2756 1942}%
\special{pa 2752 1918}%
\special{pa 2774 1906}%
\special{pa 2704 1894}%
\special{fp}%
%
\special{pn 20}%
\special{pa 2678 1970}%
\special{pa 2684 1940}%
\special{pa 2694 1908}%
\special{pa 2708 1882}%
\special{pa 2726 1856}%
\special{pa 2746 1834}%
\special{sp}%
%
\special{pn 8}%
\special{pa 2926 2008}%
\special{pa 2674 1976}%
\special{fp}%
%
\special{pn 8}%
\special{pa 2916 1980}%
\special{pa 2724 1816}%
\special{fp}%
%
\special{pn 20}%
\special{pa 3378 1784}%
\special{pa 3378 1548}%
\special{fp}%
\special{sh 1}%
\special{pa 3378 1548}%
\special{pa 3358 1616}%
\special{pa 3378 1602}%
\special{pa 3398 1616}%
\special{pa 3378 1548}%
\special{fp}%
%
\special{pn 20}%
\special{pa 3298 1552}%
\special{pa 3328 1546}%
\special{pa 3360 1544}%
\special{pa 3392 1544}%
\special{pa 3424 1544}%
\special{pa 3456 1548}%
\special{pa 3468 1550}%
\special{sp}%
%
\special{pn 8}%
\special{pa 3374 1780}%
\special{pa 3298 1546}%
\special{fp}%
%
\special{pn 8}%
\special{pa 3392 1786}%
\special{pa 3472 1554}%
\special{fp}%
%
\special{pn 20}%
\special{pa 2946 2536}%
\special{pa 2742 2668}%
\special{fp}%
\special{sh 1}%
\special{pa 2742 2668}%
\special{pa 2808 2648}%
\special{pa 2786 2638}%
\special{pa 2786 2614}%
\special{pa 2742 2668}%
\special{fp}%
%
\special{pn 20}%
\special{pa 2792 2700}%
\special{pa 2768 2680}%
\special{pa 2746 2658}%
\special{pa 2730 2632}%
\special{pa 2716 2604}%
\special{pa 2710 2582}%
\special{sp}%
%
\special{pn 8}%
\special{pa 3660 1862}%
\special{pa 4002 1360}%
\special{fp}%
%
\special{pn 20}%
\special{sh 1}%
\special{ar 2936 2008 10 10 0  6.28318530717959E+0000}%
\special{sh 1}%
\special{ar 2936 2008 10 10 0  6.28318530717959E+0000}%
%
\special{pn 20}%
\special{sh 1}%
\special{ar 2954 2536 10 10 0  6.28318530717959E+0000}%
\special{sh 1}%
\special{ar 2954 2536 10 10 0  6.28318530717959E+0000}%
%
\special{pn 20}%
\special{sh 1}%
\special{ar 3392 2728 10 10 0  6.28318530717959E+0000}%
\special{sh 1}%
\special{ar 3392 2728 10 10 0  6.28318530717959E+0000}%
%
\special{pn 20}%
\special{sh 1}%
\special{ar 3384 1780 10 10 0  6.28318530717959E+0000}%
\special{sh 1}%
\special{ar 3384 1780 10 10 0  6.28318530717959E+0000}%
%
\special{pn 20}%
\special{pa 3858 2502}%
\special{pa 4080 2602}%
\special{fp}%
\special{sh 1}%
\special{pa 4080 2602}%
\special{pa 4026 2556}%
\special{pa 4030 2580}%
\special{pa 4010 2592}%
\special{pa 4080 2602}%
\special{fp}%
%
\special{pn 20}%
\special{ar 3982 2518 142 174  0.2196967 0.9187484}%
%
\special{pn 8}%
\special{pa 3878 2500}%
\special{pa 4126 2554}%
\special{fp}%
%
\special{pn 8}%
\special{pa 3860 2510}%
\special{pa 4058 2664}%
\special{fp}%
%
\special{pn 8}%
\special{ar 4316 2564 600 266  1.5799202 2.7248086}%
%
\special{pn 8}%
\special{pa 3726 2610}%
\special{pa 3726 2610}%
\special{fp}%
%
\special{pn 20}%
\special{pa 3764 2664}%
\special{pa 3726 2618}%
\special{fp}%
\special{sh 1}%
\special{pa 3726 2618}%
\special{pa 3752 2682}%
\special{pa 3760 2660}%
\special{pa 3784 2658}%
\special{pa 3726 2618}%
\special{fp}%
%
\special{pn 20}%
\special{pa 3526 1798}%
\special{pa 3630 1596}%
\special{fp}%
\special{sh 1}%
\special{pa 3630 1596}%
\special{pa 3582 1646}%
\special{pa 3606 1644}%
\special{pa 3618 1664}%
\special{pa 3630 1596}%
\special{fp}%
%
\special{pn 20}%
\special{pa 3250 1824}%
\special{pa 3146 1596}%
\special{fp}%
\special{sh 1}%
\special{pa 3146 1596}%
\special{pa 3154 1666}%
\special{pa 3168 1646}%
\special{pa 3192 1648}%
\special{pa 3146 1596}%
\special{fp}%
\special{pa 3146 1596}%
\special{pa 3146 1596}%
\special{fp}%
%
\special{pn 20}%
\special{pa 3030 1926}%
\special{pa 2860 1742}%
\special{fp}%
\special{sh 1}%
\special{pa 2860 1742}%
\special{pa 2890 1804}%
\special{pa 2896 1782}%
\special{pa 2920 1778}%
\special{pa 2860 1742}%
\special{fp}%
%
\special{pn 20}%
\special{pa 2888 2136}%
\special{pa 2630 2116}%
\special{fp}%
\special{sh 1}%
\special{pa 2630 2116}%
\special{pa 2696 2142}%
\special{pa 2684 2120}%
\special{pa 2698 2102}%
\special{pa 2630 2116}%
\special{fp}%
%
\special{pn 20}%
\special{pa 2898 2400}%
\special{pa 2624 2436}%
\special{fp}%
\special{sh 1}%
\special{pa 2624 2436}%
\special{pa 2692 2448}%
\special{pa 2676 2430}%
\special{pa 2686 2408}%
\special{pa 2624 2436}%
\special{fp}%
%
\special{pn 8}%
\special{pa 3392 2254}%
\special{pa 3668 1844}%
\special{fp}%
%
\special{pn 8}%
\special{pa 4116 1916}%
\special{pa 4116 1916}%
\special{fp}%
%
\special{pn 8}%
\special{pa 4144 1926}%
\special{pa 4106 1926}%
\special{fp}%
\special{sh 1}%
\special{pa 4106 1926}%
\special{pa 4174 1946}%
\special{pa 4160 1926}%
\special{pa 4174 1906}%
\special{pa 4106 1926}%
\special{fp}%
%
\special{pn 8}%
\special{pa 3392 2254}%
\special{pa 2868 3056}%
\special{fp}%
%
\special{pn 20}%
\special{pa 3050 2582}%
\special{pa 2860 2800}%
\special{fp}%
\special{sh 1}%
\special{pa 2860 2800}%
\special{pa 2918 2764}%
\special{pa 2894 2760}%
\special{pa 2888 2738}%
\special{pa 2860 2800}%
\special{fp}%
%
\special{pn 20}%
\special{pa 3146 2636}%
\special{pa 2992 2884}%
\special{fp}%
\special{sh 1}%
\special{pa 2992 2884}%
\special{pa 3044 2838}%
\special{pa 3020 2838}%
\special{pa 3010 2816}%
\special{pa 2992 2884}%
\special{fp}%
%
\special{pn 20}%
\special{pa 3240 2692}%
\special{pa 3136 2946}%
\special{fp}%
\special{sh 1}%
\special{pa 3136 2946}%
\special{pa 3180 2892}%
\special{pa 3156 2898}%
\special{pa 3142 2878}%
\special{pa 3136 2946}%
\special{fp}%
%
\special{pn 20}%
\special{pa 3544 2710}%
\special{pa 3650 2938}%
\special{fp}%
\special{sh 1}%
\special{pa 3650 2938}%
\special{pa 3640 2868}%
\special{pa 3628 2890}%
\special{pa 3604 2886}%
\special{pa 3650 2938}%
\special{fp}%
\put(25.7000,-36.0000){\makebox(0,0)[lt]{($\vZ_0\,:\,$ a break point of $\rho^{-1}(b)$) }}%
%
\special{pn 8}%
\special{pa 2954 2546}%
\special{pa 2784 2710}%
\special{fp}%
%
\special{pn 8}%
\special{pa 2946 2518}%
\special{pa 2698 2582}%
\special{fp}%
\end{picture}%
\hspace{23.25truecm}}

\vspace{0.5truecm}

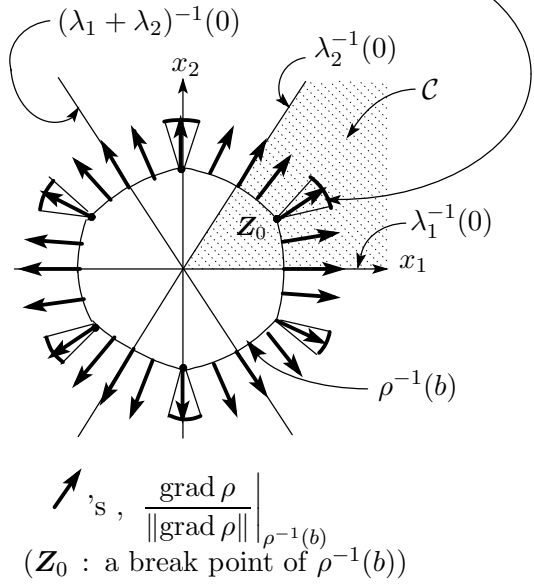
\centerline{{\bf Figure 5 $\,:\,$ The direction situation of $\displaystyle{\frac{{\rm grad}\,\rho}{\vert\vert{\rm grad}\,\rho\vert\vert}}$ (Case I)}} 


\vspace{0.65truecm}

\noindent
For any positive constant $b$, we set $\displaystyle{D_b:=\rho^{-1}((-\infty,b])}$, 
which is compact because $\rho$ is convex and proper.  
Then we can show that, for some positive constants $C_i$ ($i=1,2$), $C_1\leq\widehat F^e\leq C_2$ holds on $D_b$ 
in the weak sense because of the convexity of $\rho\vert_{\lambda_i^{-1}(0)}$ $(i=1,2$) and the definition of $\widehat F^e$.  
Also, we can show that 
$${\rm det}\left(\frac{\partial^2\rho}{\partial x_i\partial x_j}\right)=\widehat F^e$$
holds on $\mathfrak a^d$ in the weak sense.  
From these facts, we can show that 
$$C_1\leq{\rm det}\left(\frac{\partial^2\rho}{\partial x_i\partial x_j}\right)\leq C_2$$
holds on $D_b$ in the weak sense.  
On the other hand, 
$\rho\vert_{\partial D_b}$ is constant on each component of $\partial D_b$ and hence it is continuous.  
Hence, by using corollary 2 of \cite{C2}, we can show that $\rho$ is strictly convex in the classical sense on $D_{b-1}$.  
Furthermore, it follows from the arbitrariness of $b$ that $\rho$ is strictly convex in the classical sense on the whole of 
$\mathfrak a^d$.  Therefore, according to Proposition 3.6, $\rho$ is locally of class $C^{1,\alpha}$, that is, ${\rm grad}\,\rho$ is locally of class $C^{0,\alpha}$.  

Set $\rho_b:=\rho-b$.  
By operating a suitable affine transformation of $\mathfrak a^d$, 
$B^{2}(\varepsilon_b)\subset D_b\subset B^{2}(2\varepsilon_b)$ holds for some $\varepsilon_b>0$.  
Note that $\varepsilon_b\to\infty$ as $b\to\infty$.  Since $\rho_b$ satisfies 
$${\rm det}\left(\frac{\partial^2\rho_b}{\partial x_i\partial x_j}\right)=\widehat F^e\leqno{(3.26)}$$
on $D_b$ in the weak sense (that is, $\rho_b$ is the viscosity solution of $(3.26)$) 
and $\rho_b\vert_{\partial D_b}=0$, it follows from Proposition 3.4 that $\rho_b$ (hence $\rho$) is of class $C^{2,\alpha}$ ($\alpha>0$) 
on $B^{2}(\frac{\varepsilon_b}{2})$.  Since $\varepsilon_b\to\infty$ as $b\to\infty$, we see that $\rho$ is of class $C^{2,\alpha}$ on the whole of $\mathfrak a^d$.  
Furthermore, the higher regularity ($C^{\infty}$-property) of $\rho$ is shown 
by reducing this problem to the fully non-linear uniformly elliptic case (\cite{C1}) 
and using the higher regularity (which is shown by the bootstrap argument) of the $C^{2,\alpha}$-solution of fully non-linear uniformly elliptic partial differential equation 
of order two.  
Then $\omega_{\rho^h}$ is of class $C^{\infty}$ on the whole of $G^{\mathbb C}/K^{\mathbb C}$.  
From (ii) of Proposition 3.3 and $(3.25)$, for the normal complex coordinate $(z^p_1,\cdots,z^p_n)$ about $p={\rm Exp}_o(\vZ)$ 
($\vZ\in W\cdot\mathcal C$), we obtain the following expression of $(\omega_{\rho^h})_{p}^n$:
\begin{align*}
(\omega_{\rho^h})_{p}^n=&(-1)^{\frac{n(n-1)}{2}}\sqrt{-1}^n\cdot n!\cdot
{\rm det}\left(\left.\frac{\partial^2\rho^h}{\partial z^{p}_i\partial\bar z^{p}_j}\right\vert_{p}\right)\\
&\times(dz^{p}_1\wedge\cdots\wedge dz^{p}_n\wedge d\bar z^{p}_1\wedge\cdots\wedge d\bar z^{p}_n)_{p}\\
=&(-1)^{\frac{n(n-1)}{2}}\frac{\sqrt{-1}^nn!c}{4^n}\cdot(dz^{p}_1\wedge\cdots\wedge dz^{p}_n\wedge d\bar z^{p}_1\wedge\cdots\wedge d\bar z^{p}_n)_{p}.
\end{align*}
Also, $\Omega_p$ is expressed as 
$$\Omega_{p}=(dz^{p}_1\wedge\cdots\wedge dz^{p}_n)_{p},$$
Hence we obtain 
$$(\omega_{\rho^h})_{p}^n=(-1)^{\frac{n(n-1)}{2}}\frac{\sqrt{-1}^nn!c}{4^n}\cdot\Omega_{p}\wedge\overline{\Omega}_{p}.\leqno{(3.27)}$$
Since this relation $(3.27)$ holds at any point $p$ of ${\rm Exp}_o(W\cdot\mathcal C)$, it follows from the $G$-invariance of $\omega_{\rho^h}$ and $\Omega$ that 
$(3.27)$ holds for any point $p$ of $G\cdot{\rm Exp}_o(W\cdot\mathcal C)(=(G^{\mathbb C}/K^{\mathbb C})_{\rm reg})$.  
Furthermore, since $\omega_{\rho^h}^n$ and $\Omega\wedge\overline{\Omega}$ are continuous on $G^{\mathbb C}/K^{\mathbb C}$, $(3.27)$ holds on the whole of 
$G^{\mathbb C}/K^{\mathbb C}$.  
If $c=2^n$, then we have 
$$(\omega_{\rho^h})_{p}^n=(-1)^{\frac{n(n-1)}{2}}\left(\frac{\sqrt{-1}}{2}\right)^n\cdot n!\cdot
\Omega_{p}\wedge\overline{\Omega}_{p}\quad\,\,(p\in G^{\mathbb C}/K^{\mathbb C}).\leqno{(3.28)}$$
Thus the quadruple $({\bf J},\beta_{\rho^h},\omega_{\rho^h},\Omega)$ gives a $C^{\infty}$-Calabi-Yau structure on \newline
$G^{\mathbb C}/K^{\mathbb C}$.  \qed


\vspace{1truecm}

{\small 
\rightline{Department of Mathematics, Faculty of Science}
\rightline{Tokyo University of Science, 1-3 Kagurazaka}
\rightline{Shinjuku-ku, Tokyo 162-8601 Japan}
\rightline{(koike@rs.tus.ac.jp)}
}

\end{document}